\documentclass[11pt]{article}

\usepackage{latexsym,amsmath,amsfonts,amssymb,amstext,
euscript,array,enumerate}
\usepackage{color}
\usepackage{comment}

\usepackage{amsthm}
\usepackage{bm}
\allowdisplaybreaks
\usepackage[colorlinks=true,linkcolor=blue,urlcolor=blue]{hyperref}

\numberwithin{equation}{section}

\theoremstyle{plain}
\newtheorem{theorem}{Theorem}[section]
\newtheorem{lemma}[theorem]{Lemma}
\newtheorem{proposition}[theorem]{Proposition}

\newtheorem{Definition}{Definition}[section]

\newtheorem{Remark}{Remark}[section]

\numberwithin{equation}{section}

\theoremstyle{Remark}

\def\esssup_#1{\underset{#1}{\mathrm{ess\,sup\, }}}
\def\essinf_#1{\underset{#1}{\mathrm{ess\,inf\, }}}

\def\qed{{\hfill\hbox{\enspace${ \square}$}} \smallskip}
\def\sqr#1#2{{\vcenter{\vbox{\hrule height .#2pt \hbox{\vrule
 width .#2pt height#1pt \kern#1pt \vrule
width .#2pt} \hrule height .#2pt}}}}
\def\square{\mathchoice\sqr54\sqr54\sqr{4.1}3\sqr{3.5}3}

\def\ds{\begin{displaystyle}}
\def\eds{\end{displaystyle}}

\def\<{\langle }
\def\>{\rangle }

\def \N{\mathbb{N}}
\def \R{\mathbb{R}}

\def \E{\mathbb{E}}
\def \F{\mathbb{F}}

\def \P{\mathbb{P}}

\def \Vc{{\cal V}}

\def \Vc{{\cal V}}

\def \eps{\varepsilon}

\def\beqs{\begin{eqnarray*}}
\def\enqs{\end{eqnarray*}}
\def\beq{\begin{eqnarray}}
\def\enq{\end{eqnarray}}

\newcommand{\spernuxa}[1]{\mathbb{E}^{x,a_0,a_\Gamma}_{\bm\nu} \left[ #1 \right]}

\newcommand{\spernuxaepsilon}[1]{\mathbb{E}^{x,a_0,a_\Gamma}_{\bm\nu^{\varepsilon}} \left[ #1 \right]}                          
\newcommand{\sperxa}[1]{\mathbb{E}^{x,a_0,a_\Gamma} \left[ #1 \right]}                               

\DeclareMathAlphabet{\mathonebb}{U}{bbold}{m}{n}                           %
\newcommand{\one}{\ensuremath{\mathonebb{1}}}                               

\allowdisplaybreaks

\title{Constrained BSDEs driven by a non quasi-left-continuous random measure and  optimal control of PDMPs on bounded domains}

\date{}
\author{Elena BANDINI\thanks{Universit\`{a} degli Studi di Milano-Bicocca, Via Roberto Cozzi, 55, 20125 Milano  Italy; e-mail: \texttt{elena.bandini@unimib.it}}
}

\begin{document}

\maketitle

\begin{abstract}
We consider an  
 optimal control problem for  piecewise deterministic Markov processes (PDMPs) on a bounded state space. 
The control problem under study is very general: a
  pair of  controls acts continuously on the deterministic flow and on  the two transition measures (in the interior and from the boundary of the domain) describing  the jump dynamics of the process.
 For this class of control problems, the value function can  be characterized as the unique viscosity solution to the corresponding fully-nonlinear Hamilton-Jacobi-Bellman equation with a non-local type boundary condition.
 
 By means of the recent control  randomization method, we are able to  provide a probabilistic representation for the value function in terms of a   constrained backward stochastic differential equation (BSDE), known as nonlinear Feynman-Kac formula. 
 This result considerably extends the existing literature, 
 where only  the  case with no jumps from the   boundary is considered. 
 The  additional 
 boundary jump mechanism is described in terms of a non quasi-left-continuous random measure and induces  
predictable jumps 
in the PDMP's dynamics.
 The existence and uniqueness results for BSDEs  driven by 
such a random measure  are non trivial, even in the unconstrained case, as emphasized in the recent work \cite{BandiniBSDE}.

\end{abstract}


\medskip


 \section{Introduction}\label{Sec_introduction}
 In this paper we prove that the value function of an infinite-horizon  optimal control problem for  piecewise deterministic Markov processes on bounded domains can be represented in terms of a suitable backward stochastic differential equation. 
 Piecewise deterministic Markov processes,  introduced in   \cite{Da}, evolve
 by means of random jumps at random times,
 while the behavior between jumps
 is described  by a deterministic flow. 
  We consider   optimal control problems of PDMPs where  the control acts continuously on the jump dynamics as well as
 on the deterministic flow.
We deal with PDMPs with  bounded state space: whenever the process hits the boundary, it immediately jumps into the interior of the domain.
 Control  problems for this type of processes arise in many contexts, among which operations research, engineering systems and management science, see \cite{Da} for a detailed overview.
Our aim is to represent the value function  by means of an appropriate BSDE.
It is worth mentioning that the	probability measures   describing the distribution of the controlled PDMP are in general not absolutely continuous with respect to the law of a given, uncontrolled process (roughly speaking, the control problem is non-dominated). 
This is  reflected 
 	  in the fully nonlinear character of the associated HJB equation, and prevents the use of standard BSDE techniques.
For this reason, we shall extend to the present framework the so-called \emph{randomization method},  recently introduced by   \cite{KhPh} 
in the diffusive context,
to represent  the solutions of 
 fully nonlinear integro-partial differential equations by means of a new class of BSDEs with nonpositive jumps.
The extension of the randomization approach
 to our PDMPs optimal control problem  is particularly delicate due to the presence of the  jump mechanism from the boundary. 
 Indeed, since the  jumps  from the boundary happen at predictable times, the associated BSDE  turns out to be driven by a non quasi-left-continuous random measure.
  For such general backward equations, 
  the existence and uniqueness of a solution  is particularly tricky, and 
counterexamples  can be obtained even in simple  cases, see  \cite{CFJ}.
 

 Let us  describe our setting in more detail.
Let $E$ be an open bounded subset of $\R^d$, with Borel $\sigma$-algebra $\mathcal E$. The set $E$ corresponds to the PDMP state space.  
 Roughly speaking, a controlled PDMP on 
 $(E,\mathcal E)$  is described by specifying  its  local characteristics, namely  a 
vector field   $h(x,a_0)$, a jump rate $\lambda(x,a_0)$, and two transition probability measures $Q(x,a_0,dy)$
 and  $R(x,a_{\Gamma},dy)$ prescribing the positions of the process at the jump times, respectively starting from the interior and from the boundary of the domain.
The local characteristics   depend on some initial value $x \in E$ and on  the parameters  $a_0\in A_0$,  $a_{\Gamma}\in A_{\Gamma}$, where $(A_0, \mathcal{A}_0)$ and $(A_{\Gamma}, \mathcal{A}_{\Gamma})$ are two general measurable spaces,  denoting respectively the space  of control actions in the interior and on  the boundary  of the domain.
 The control procedure consists in  
 choosing a  pair of  strategies: 
a \emph{piecewise open-loop policy} controlling  the motion in the interior of the domain, i.e. a  measurable function only depending on the last jump time $T_n$ and post jump position $E_n$, and a boundary control belonging  to the set of \emph{feedback policies}, that only depends on the position of the process just before the jump time. The above formulation of the control problem is used in many papers as well as books, see for instance  
 \cite{CoDu},  \cite{Da}. 
 The class of admissible control laws $\mathcal A_{ad}$ will be  the set of all $\mathcal A_0 \otimes \mathcal A_\Gamma$-measurable maps $\bm \alpha =(\alpha^0, \alpha^\Gamma)$, with  $\alpha^{\Gamma}: \partial E \rightarrow A_{\Gamma}$, and $\alpha^0: [0,\,\infty) \times E \rightarrow A_0$ such that 
 \begin{eqnarray*}
  \alpha^0_t=\alpha^0_0(t,x)\,\one_{[0,\,T_1)}(t)+ \sum_{n=1}^{\infty}\alpha^0_n(t-T_n,E_{n})\,\one_{[T_n,\,T_{n+1})}(t).
 \end{eqnarray*}  
The controlled process $X$ is defined as 
  \begin{equation*}
  X_t=
  \left\{
  \begin{array}{ll}
    \phi^{\alpha^0}(t, x) \quad &\textup{if}\,\,t \in [0,\, T_{1}),\\
  \phi^{\alpha^0}(t- T_n, E_n)\quad &\textup{if}\,\, t \in [T_n,\, T_{n+ 1}),\,\, n \in \N \setminus \{0\},
  \end{array}	
  \right.
  \end{equation*}
 where $\phi^{\alpha^0}(t, x)= \phi(t,x,\alpha^0_t)$
  is the unique solution to the ordinary differential equation on $\R^d$
  $$
  \dot x(t) = h(x(t), \alpha^0(t)), \quad x(0)=x.
  $$
For every starting point $x \in E$ and for each $\bm \alpha \in \mathcal A_{ad}$, 
one can introduce the unique probability measure $\P^{x}_{\bm \alpha}$ 
such that
 the conditional survival function of
 the   inter-jump  times
 and the distribution of the
 post jump positions of $X$ under $\P^{x}_{\bm \alpha}$ 
 are given by \eqref{dyn1}-\eqref{dyn2}-\eqref{dyn3}.
We denote by  $\E^{x}_{\bm\alpha}$ the expectation under $\P^{x}_{\bm\alpha}$.
 In the classical infinite-horizon control problem the goal is to minimize over all control laws $\bm\alpha$ a functional cost of the form
\[
 J(x,\bm\alpha) = \E^{x}_{\bm\alpha}\left[\int_{(0,\infty)} e^{-\delta\,s}\, f(X_{s},\alpha^0_s)\,ds+ \int_{(0,\infty)} e^{-\delta\,s}\, c(X_{s-},\alpha^{\Gamma}(X_{s-}))\,d p^\ast_s\right],
\]
where $f$  is a given real function on  $\bar E \times A_0$ representing the running cost, $c$  is a given real function on  $\partial E \times A_{\Gamma}$ that provides a cost every time the process hits the boundary,  $\delta \in (0,\,\infty)$ is a  discount factor, 
while the process $p^\ast_s$
counts the number of times the boundary is hit (see \eqref{pstar}).
 The value function of the control problem is defined in the usual way:
 \begin{equation}\label{Sec:PDP_value_functionINTRO}
 V(x) = \inf_{\bm\alpha \in \mathcal{A}_{ad}}J(x,\bm\alpha),\quad x \in E.
 \end{equation} 
 Under suitable assumptions on the cost functions $f,c$, and on the local characteristics $h, \lambda, Q, R$, 
 $V$ is known to be  the unique continuous viscosity solution on $[0,\,\infty) \times \bar E$ of the Hamilton-Jacobi-Bellman (HJB) equation with boundary non-local condition: 
 \begin{eqnarray}\label{Intro_HJB}
  \left\{
  \begin{array}{ll}
 \delta v(x)= \inf_{a_0 \in A_0} \left(h(x,a_0)\cdot \nabla v(x) + \lambda(x,a_0) \int_E (v(y)-v(x))\,Q(x,a_0,dy)+ f(x,a_0)\right),\,\, x \in E,\label{Sec:PDP_HJBINTRO}\\
 v(x)= \min_{a_{\Gamma} \in A_{\Gamma}}(\int_E (v(y)-v(x))\,R(x,a_{\Gamma},dy)+c(x,a_\Gamma)), \,\, x \in\partial E.
  \end{array}
  \right.
 \end{eqnarray}
Our aim is to represent the value function $V$ by means of an appropriate BSDE.
We are interested in the general case where 
the	probability measures   $\{\P^x_{\alpha}\}_{\alpha}$ describing the distribution of the controlled process are  not absolutely continuous with respect to the law of a given, uncontrolled process.
 Probabilistic formulae for the value function for non-dominated models have been discovered only in recent years.
 In this sense, 
a key role is played by  the randomization method, firstly introduced 
in
 \cite{KhPh} 
 to
 represent the solutions of 
 fully nonlinear integro-partial differential equations related to the classical optimal control for diffusions, 
 and later extended to other types of control problems, 
   see for instance \cite{KhMaPhZh}, \cite{ElKh14a},  \cite{Co-Chou}, 
  \cite{CoFuPh}, 
   \cite{BaCoFuPh1}. 
In the non-diffusive  framework, 
  the correct formulation of the randomization method  requires some efforts and different techniques from 
  the diffusive case, since the controlled process  is naturally described only in terms of its local characteristics and not as a solution to some stochastic differential equation. 
 A  first step in  the generalization of the randomization method to the non-diffusive framework was done in \cite{BandiniFuhrman},
  where a probabilistic representation for the value function associated to an optimal control problem for pure jump Markov processes was provided;
afterwards, the randomization techniques 
 have been implemented  in \cite{BandiniPDMPsNoBordo} to solve PDMPs optimal control problems on unbounded state spaces  (notice that  in both \cite{BandiniFuhrman} and  \cite{BandiniPDMPsNoBordo} the jump measure of the controlled state process is quasi-left continuous).  
 In the present paper we are interested to extend those results to the case of optimal control problems for PDMPs on bounded state spaces,  
 where additional
 forced jumps appear whenever the process hits the boundary.
 The   jump mechanism from the boundary plays a fundamental role as it 
 leads, among other things,   to the study  of  
  BSDEs   
  driven by a non quasi-left-continuous random measure.
Only 
recently,  
some results have been obtained on this subject,  see \cite{CohenElliottPearce}, \cite{CohenElliott}, \cite{BandiniBSDE};  
 in particular, in \cite{BandiniBSDE}  well-posedness  is proved for unconstrained BSDEs in a general non-diffusive framework, under  a 
 specific condition involving the Lipschitz constants of the BSDE generator  and the size of the predictable jumps. In the present paper we extend the results in \cite{BandiniBSDE} 
 to our class of constrained BSDEs.

We now describe   the randomization  approach in our framework. 
 The fundamental idea  consists in 
the so-called \emph{randomization of the control}: 
roughly speaking, we  replace the state trajectory and the associated pair of controls $(X_s,\alpha^0_s, \alpha^{\Gamma}_s)$ by an (uncontrolled) PDMP $(X_s,I_s, J_s)$.   The process  $I$ (resp. $J$) is chosen to be a pure  jump process with values in the space of control actions $A_0$ (resp. $A_{\Gamma}$), with an intensity $\lambda_0(db)$ (resp. $\lambda_{\Gamma}(dc)$),  which is arbitrary but finite and with full support.
 In particular, the PDMP $(X,I, J)$ is  constructed on a new probability space by means of a different triplet of local characteristics and takes values on the enlarged space $E \times A_0 \times A_{\Gamma}$ (or, equivalently, by assigning the compensator $\tilde p(ds\,dy\,db\,dc)$). For any  starting point $(x,a_0,a_{\Gamma})$  in $E \times A_0 \times A_{\Gamma}$, we  denote by $\P^{x,a_0, a_{\Gamma}}$ the corresponding law.
 At this point  we introduce  an auxiliary optimal control problem where we control the intensity of the processes $I$ and $J$: using a Girsanov's type theorem for point processes, for any pair of predictable, bounded and positive processes $(\nu^0_t(b), \nu^\Gamma_t(c))$,  we construct a probability measure $\P^{x,a_0,a_{\Gamma}}_{\nu^0,\nu^\Gamma}$ under which the compensator of $I$ (resp. $J$) is given by $\nu^0_t(db)\,\lambda_0(db)\,dt$ (resp. $\nu^\Gamma_t(dc)\,\lambda_{\Gamma}(dc)\,dt$). 
 It is worth mentioning that
  the applicability of the Girsanov theorem  to the present framework, i.e. when the compensator $\tilde p$  is  a non quasi-left-continuous random measure, requires the validity of  an additional condition involving the intensity control fields  $(\nu^0,\nu^\Gamma)$ and  the predictable jumps of $\tilde p$, see  \eqref{cond_T4.5_J}. 
The correct formulation of the randomized control problem has to take into account this  latter constraint.

 The aim of the new control problem (called \emph{randomized} or \emph{dual} control problem) is to minimize the functional
 \begin{equation}\label{Sec:PDP_functional_dualINTRO}
 J(x,a_0,a_{\Gamma}, \nu^0,\nu^\Gamma) = \E^{x,a_0,a_{\Gamma}}_{\nu^0,\nu^\Gamma}\left[\int_{(0,\infty)} e^{-\delta\,s}\, f(X_{s},I_s)\,ds+\int_{(0,\infty)} e^{-\delta\,s}\, c(X_{s-},J_{s-})\,d p^\ast_s\right]
 \end{equation}
 over all possible choices of $\nu^0,\nu^\Gamma$. 
 Firstly, we  give a
  probabilistic  representation of the value function of the randomized control problem, denoted $V^\ast(x,a_0,a_{\Gamma})$, in terms of of a well-posed constrained BSDE. This latter  is an equation over  infinite horizon of  the form \eqref{Sec:PDP_BSDE}
with the sign constraints \eqref{Sec:PDP_BSDE_constraint1}-\eqref{Sec:PDP_BSDE_constraint2}.
The random measure $q = p-\tilde p$ driving the BSDE is the compensated  measure associated to the jumps of $(X,I, J)$. In particular,  the compensator $\tilde p$ has predictable jumps  $\tilde p(\{t\}\times dy\,db\,dc)= \one_{X_{t-} \in \partial E}$.
Equation \eqref{Sec:PDP_BSDE}-\eqref{Sec:PDP_BSDE_constraint1}-\eqref{Sec:PDP_BSDE_constraint2}  is thus driven by a non quasi-left-continuous random measure; the associated   well-posedness results are obtained by means of a penalization approach,
by suitably extending the  recent existence and uniqueness  theorem obtained in \cite{BandiniBSDE} for unconstrained BSDEs.
Once we achieve the existence and uniqueness of a maximal solution to \eqref{Sec:PDP_BSDE}-\eqref{Sec:PDP_BSDE_constraint1}-\eqref{Sec:PDP_BSDE_constraint2}, we prove that its component $Y^{x,a_0,a_{\Gamma}}$ at the initial time represents the randomized value function, i.e. 
$ Y^{x,a_0,a_{\Gamma}}_0= V^{\ast}(x,a_0,a_{\Gamma})$.
All this is collected in Theorem \ref{Sec:PDP_Thm_ex_uniq_maximal_BSDE}.
Then, we aim at proving that  $ Y^{x,a_0,a_{\Gamma}}_0$
also  provides a nonlinear Feynman-Kac representation to the value function \eqref{Sec:PDP_value_functionINTRO} of our original  optimal control problem.
 To this end, we introduce the deterministic real function on $E \times A_0\times A_{\Gamma}$ defined by 
$v(x,a_0,a_{\Gamma}):=Y_0^{x,a_0,a_{\Gamma}}$. In Theorem \ref{Sec:PDP_THm_Feynman_Kac_HJB} we prove that $v$ does not depend on its two last arguments, is a bounded and continuous function on $E$, and   that $v(X_s)= Y_s^{x,a_0,a_\Gamma}$ for all $s \geq 0$.
Then, we show that $v$ is a viscosity solution to \eqref{Sec:PDP_HJBINTRO}, so that,  
 by the uniqueness of the solution to the HJB equation \eqref{Sec:PDP_HJBINTRO},  we can conclude that 
 \begin{equation}\label{Sec:PDP_rapprINTRO}
Y_0^{x,a_0,a_{\Gamma}}=V^\ast(x,a_0,a_{\Gamma})=V(x).
 \end{equation}
 This   constitutes the main result of the paper  and is stated in Theorem \ref{Sec:PDP_THm_Feynman_Kac_HJB_2}.
Formula \eqref{Sec:PDP_rapprINTRO} gives  the desired BSDE representation of the value function for the original control problem. 
 This nonlinear Feynman-Kac formula can be used to design
 algorithms
 based on  the numerical approximation of the solution to the constrained BSDE \eqref{Sec:PDP_BSDE}-\eqref{Sec:PDP_BSDE_constraint1}-\eqref{Sec:PDP_BSDE_constraint2}, and therefore to get
 probabilistic numerical approximations  for the value function of the considered    optimal control problem.
 Recently, 
 numerical schemes for constrained BSDEs have been  proposed and analyzed in the diffusive framework, see 
  \cite{KhLaPha}, and  
in the PDMPs context  as well, see \cite{AbergelHurePham}.
 \normalcolor
 
 The paper is organized as follows. 
  In Section 
 \ref{Sec:PDP_Sec_control_problem} 
 we introduce the optimal control \eqref{Sec:PDP_value_functionINTRO}, and we discuss its solvability.
In Section \ref{Sec:PDP_Section_dual_control}
 we formulate   the randomized optimal control problem  \eqref{Sec:PDP_functional_dualINTRO}.
 In Section \ref{Sec:PDP_Sec_ConstrainedBSDE} we introduce the constrained BSDE \eqref{Sec:PDP_BSDE}-\eqref{Sec:PDP_BSDE_constraint1}-\eqref{Sec:PDP_BSDE_constraint2} over infinite horizon,
 we show that   it admits a unique maximal solution $(Y,Z,K)$ in a certain class of processes, and that  $Y_0$  coincides with the value function of the randomized optimal control problem. 
Then, in Section \ref{Sec:PDP_Section_nonlinear_IPDE} we prove that 
$Y_0$ 
also  provides a viscosity solution to \eqref{Sec:PDP_HJBINTRO}. 
Finally, some  technical results 
are collected in the Appendix.

\section{Optimal  control of PDMPs on bounded domains}\label{Sec:PDP_Sec_control_problem}

In this  section we formulate
the  optimal control problem for  piecewise deterministic Markov processes on bounded domains, and we discuss its solvability.
 The PDMP state space $E$ is an open  bounded subset of  $\R^d$, and $\mathcal E$ the corresponding Borel $\sigma$-algebra. Moreover, we introduce two 
 Borel spaces (i.e. 	topological spaces homeomorphic to  Borel subsets of  compact metric spaces) $A_0$, $A_\Gamma$, endowed with their Borel $\sigma$-algebras $\mathcal{A}_0$ and $\mathcal{A}_\Gamma$, that are respectively the space of control actions in the interior and on the boundary of the domain. 
 Given a topological space $F$, in the sequel we will denote by  $\mathbb C_b(F)$ (resp. $\mathbb C^1_b(F)$) the set of all bounded continuous functions (resp.  all bounded  differentiable functions whose derivative is continuous) on $F$. 

A controlled PDMP on 
 $(E,\mathcal E)$  is described by means of a set of local characteristics $(h, \lambda,Q, R)$, with $h, \lambda$  functions on $\bar E \times A_0$, and $Q$, $R$  probability transition measures in $E$  respectively from $\bar E \times A_0$ and from $\partial E \times A_\Gamma$.  We assume the following.

 \medskip
 
 \noindent \textbf{(H$\textup{h$\lambda$QR}$)}
 \begin{itemize}
 	\item[(i)]
 	$h  : \bar E \times A_0 \rightarrow \R^d$ and $\lambda: \bar E \times A_0 \rightarrow \R_+$ are   
 	continuous and bounded functions, Lipschitz continuous on $\bar E$, uniformly in $A_0$.	
 	\item[(ii)]
 	$Q$ (resp . $R$)  maps 
 $\bar E \times A_0$ 
 	(resp. $\partial E \times A_\Gamma$)
 	into the set of probability measures on $(E, \mathcal E)$, and is a 
 	continuous
 	stochastic  kernel.
 	Moreover, for all $v\in\mathbb C_b(E)$, the maps  $(x, a_0) \mapsto \int_{E}v(y)\,Q(x,a_0,dy)$  and   $(x, a_\Gamma)\mapsto \int_{E}v(y)\,R(x,a_\Gamma,dy)$  are Lipschitz continuous in $x$, uniformly in $a_0 \in A_0$ and  in $a_\Gamma \in A_\Gamma$, respectively.
 \end{itemize}
    We construct the process $X$ on $E$ in a canonical way.
To this end, let $\Omega'$ be the set of sequences $\omega'=(t_n, e_n)_{n \geq 1}$ contained in $((0,\,\infty) \times E\cup \{ (\infty, \Delta)\}$, where $\Delta \notin E$, is  adjoined to $E$ as an isolated point, such that $t_n \leq t_{n+1}$, and $t_n < t_{n+1}$ if $t_n < \infty$.
 We set $\Omega = E\times \Omega'$, where $\omega =(x,\omega') = (x,t_1, e_1, t_2, e_2,  ...)$.
On the sample space $\Omega$  we define the canonical functions $T_n : \Omega \rightarrow (0,\,\infty]$, $E_n : \Omega \rightarrow E \cup \{\Delta\}$ as follows: 
 $T_0(\omega)=0$, $E_0(\omega)=x$, and  for 
 $n \geq 1$,
 $ T_n (\omega )= t_n$,  
  $E_n(\omega)=e_n$, and $T_\infty  (\omega )= \lim_{n\to\infty} t_n$.
We also introduce, for any $B \in \mathcal E$,  the counting  process $N(s,B)= \sum_{n \in \N} \one_{T_n \leq s}\one_{E_n \in B}$ and the associated integer-valued random measure    on $(0,\,\infty) \times E$ 
$$
 p(ds\,dy)= \sum_{n\in \N} \one_{\{T_n,E_n\}}(ds\,dy).
$$
The class of admissible control maps $\mathcal A_{ad}$ is the set of all $\mathcal A_0 \otimes \mathcal A_\Gamma$-measurable maps $\bm \alpha =(\alpha^0, \alpha^\Gamma)$, where  $\alpha^0: [0,\,\infty) \times E \rightarrow A_0$ is a piecewise open-loop  function of the form
 \begin{align*}
  \alpha^0_t &=\alpha^0_0(t,x)\,\one_{[0,\,T_1)}(t)+ \sum_{n=1}^{\infty}\alpha^0_n(t-T_n,E_{n})\,\one_{[T_n,\,T_{n+1})}(t)
 \end{align*}
and  $\alpha^{\Gamma}:  \partial E \rightarrow A_{\Gamma}$ is a feedback policy.
 We define  the controlled  process $X: \Omega \times [0,\,\infty) \rightarrow \bar E \cup \{\Delta\}$ setting
   \begin{equation}\label{Sec:PDP_controlledX}
   X_t=
   \left\{
   \begin{array}{ll}
   \phi^{\alpha_0^0}(t, x) \quad &\textup{if}\,\,t \in [0,\, T_{1}),\\
   \phi^{\alpha_n^0}(t- T_n, E_n)\quad &\textup{if}\,\, t \in [T_n,\, T_{n+ 1}),\,\, n \in \N \setminus \{0\},
   \end{array}	
   \right.
   \end{equation}
 where $\phi^U(t, x)$, with $U$ any $\mathcal A_0$-measurable function, 
 is the unique solution to the ordinary differential equation
 $\dot y(t) = h(y(t), U(t))$, $y(0)=x \in E$. 
Finally, we introduce  the process 
\begin{equation}\label{pstar}
p^\ast_s:= \sum_{n=1}^{\infty}\one_{\{s \geqslant T_n\}}\,\one_{\{X_{T_{n}-} \in \partial E\}},
\end{equation}
that counts the number of times that the process hits the boundary.

Set $\mathcal F_0 = \mathcal E \otimes \{\emptyset, \Omega'\}$
 and, for all $t \geq 0$,  $\mathcal{G}_t= \sigma(p((0,s] \times B): s \in (0,t], B \in \mathcal E)$. For all $t$, let $\mathcal F_t$ be the $\sigma$-algebra generated by $\mathcal F_0$ and $\mathcal G_t$.
In the following all the concepts of measurability for stochastic processes will refer to the right-continuous, natural  filtration $\mathbb F = (\mathcal F_t)_{t \geq 0}$.  By the symbol $\mathcal P$ we will denote the $\sigma$ algebra of $\mathbb F$-predictable subsets of $[0,\,\infty) \times \Omega$.

For every starting point $x \in E$ and for each $\bm\alpha \in \mathcal A_{ad}$, by Theorem 3.6 in \cite{J},
there exists a unique probability measure on $(\Omega, \mathcal F_{\infty})$,  denoted by $\P^{x}_{\bm \alpha}$, such that its restriction to $\mathcal F_0$ is  the Dirac measure concentrated at  $x$,  and
 the  $\mathbb F$-compensator 
 under $\P^{x}_{\bm\alpha}$ of the measure $p(ds\,dy)$ is
 \begin{align*}
 \tilde{p}^{\bm\alpha}(ds\,dy)&=
 \sum_{n=1}^{\infty}\one_{[T_n,\,T_{n+1})}(s)\,\one_{\{\phi^{\alpha_n^0}(s-T_n,  E_n)\in E\}}\\
 & \qquad  \cdot \lambda(\phi^{\alpha_n^0}(s-T_n,  E_n)),\alpha_n^0(s-T_n,  E_n)
 )\,Q(\phi^{\alpha_n^0}(s-T_n,  E_n)),\alpha_n^0(s-T_n,  E_n)
 , dy)\,ds\\
 &+ \sum_{n=1}^{\infty}\one_{[T_n,\,T_{n+1})}(s)\,\one_{\{\phi^{\alpha_n^0}(s-T_n,  E_n)\in \partial E\}}\\
 & \qquad  \cdot R(\phi^{\alpha_n^0}(s-T_n,  E_n)),\alpha^\Gamma(\phi^{\alpha_n^0}(s-T_n,  E_n)), dy)\,d p ^\ast_s.
 \end{align*}
 Arguing as in Proposition 2.2 in \cite{BandiniPDMPsNoBordo}, one can easily see that under $\P^{x}_{\bm \alpha}$ the process $X$ in \eqref{Sec:PDP_controlledX} is 
 Markovian  with respect to  $\mathbb F$.
 In particular, for every $n \geq1$,
 the conditional survival function of
 the   inter-jump  time $T_{n+1}-T_n$ on $\{T_n < \infty\}$
is
 \begin{align}\label{dyn1}
 &\P^x_{\bm\alpha}(T_{n+1}-T_n > s\,|\,\mathcal F_{T_n})=
 \exp\left(-\int_{T_n}^{T_n +s} \lambda (\phi^{\alpha^0}(r-T_n,X_{T_n}),\alpha^0_n(r-T_n,X_{T_n}))\,dr\right)\,
\one_{\phi^{\alpha^0}(s,X_{T_n}) \in E},
  \end{align}
  and the distribution of the
 post jump position $X_{T_{n+1}}$  on
 $\{T_n < \infty\}$ are
 \begin{align}
 &\P^x_{\bm\alpha}(X_{T_{n+1}} \in B 
 |\,\mathcal F_{T_n},\,T_{n+1}, \phi^{\alpha^0}(T_{n+1}-T_n,X_{T_n}) \in  E)\notag\\
 &=Q(\phi^{\alpha^0}(T_{n+1}-T_n,X_{T_n}),\alpha^0_n(T_{n+1}-T_n,X_{T_n}), B),\quad \forall B \in \mathcal E,\label{dyn2}\\
 &\P^x_{\bm\alpha}(X_{T_{n+1}} \in B 
 |\,\mathcal F_{T_n},\,T_{n+1},\phi^{\alpha^0}(T_{n+1}-T_n,X_{T_n}) \in \partial E)\notag\\
 &=R(\phi^{\alpha^0}(T_{n+1}-T_n,X_{T_n}),\alpha^{\Gamma}(\phi^{\alpha^0}(T_{n+1}-T_n,X_{T_n})), B),\quad \forall B \in \mathcal E.\label{dyn3}
 \end{align}
The infinite horizon control problem consists in minimizing over all control laws $\bm\alpha$ a cost functional of the following form:
\[
 J(x,\bm\alpha) = \E^{x}_{\bm\alpha}\left[\int_{(0,\infty)} e^{-\delta\,s}\, f(X_{s},\alpha^0_s)\,ds+ \int_{(0,\infty)} e^{-\delta\,s}\, c(X_{s-},\alpha^{\Gamma}(X_{s-}))\,d p^\ast_s\right],
\]
where $f$  is a given real function on  $\bar E \times A_0$ representing the running cost, $c$  is a given real function on  $\partial E \times A_{\Gamma}$ that associates a cost to hitting the active boundary,  $\delta \in (0,\,\infty)$  is a  discounting factor.
  The value function of the control problem is defined in the usual way:
 \begin{equation}\label{Sec:PDP_value_function}
 V(x) = \inf_{\bm\alpha \in \mathcal{A}_{ad}}J(x,\bm\alpha),\quad x \in E.
 \end{equation}
We ask that $f$ and $c$ verify the following conditions.

 \medskip

 \noindent \textbf{(H$\textup{fc}$)}\quad
 $f:  \bar E \times A_0 \rightarrow \R_+$ (resp. $c: \partial E \times A_\Gamma \rightarrow \R_+$) 
is a   continuous and bounded  function,  Lipschitz  continuous on $\bar E$ (resp. on $\partial E$), uniformly in $A_0$ (resp. $A_\Gamma$). 
In particular,
 	\begin{equation*}
 	\left\{
 	\begin{array}{ll}
 	|f(x,a)| \leqslant M_f,\quad \forall x\in  \bar E,\,\,a \in A_0, \\
 	|c(x,a)| \leqslant M_c,\quad \forall x \in \partial E,\,\,a \in A_\Gamma.
 	\end{array}
 	\right. 
 	\end{equation*}
Moreover,    set
$ t^{\alpha^0}_\ast(x) := \inf\{ t \geqslant 0 : \phi^{\alpha^0}(t,x) \in \partial E,\, x \in E \}$, and 
$E_{\varepsilon} := \left\{ x \in E: \inf_{\alpha^0 \in \mathcal{A}_0} t_\ast^{\alpha^0}(x) \geqslant \varepsilon \right\}$. 
We will consider the following assumption. 

\medskip 

\noindent \textbf{(H$0$)}  
    There exists $\varepsilon >0$  such that $R(x,\alpha, E_\varepsilon)=1$ for all $x \in \partial E$ and $\alpha \in \mathcal A_\Gamma$.
\begin{Remark}\label{Rem_Nt}
Roughly speaking, condition \textup{\textbf{(H$0$)}} says that the state process always jumps from the boundary to points of the interior of the domain whose distance from the boundary (as measured by the boundary hitting time $t^{\alpha^0}_\ast$) are uniformly bounded away from zero.
Together with the boundedness assumption of the jump rate $\lambda$ in \textup{\textbf{(H$\textup{h$\lambda$QR}$)}-(i)}, it insures 
that, for every starting point $x \in E$ and 
admissible control $\bm\alpha\in \mathcal A_{ad}$, we have (see the proof of Proposition 24.6 in \cite{Da})
\begin{align}
\E^{x}_{\bm\alpha}\left[p^{\ast}_t\right] &\leqslant  \frac{t}{\varepsilon} + 1 =: C^\ast(t), \quad \forall \,\,t\geq0.\label{conv_exp_sper_pstar}
\end{align} 
By the integration by parts formula for processes of finite variation (see e.g., Proposition 4.5 in \cite{Re-Yor}), this implies in particular that
\begin{equation}\label{Cdelta}
\E^{x}_{\bm\alpha}\left[\int_{(0,\infty)} e^{-\delta\,t}\,dp^\ast_t\right]
\leqslant \frac{1}{\eps} + 1 =: C^\ast.
\end{equation}
\end{Remark}
Finally, we impose the following standard non-degeneracy assumptions, that allow to avoid difficulties arising from trajectories tangent to the boundary, see \cite{Barles94} for more details.

\medskip

\noindent \textbf{(HBB)} \quad For all $x \in \partial E$, if there exists $a_0 \in A_0$ such that $-h(x,a_0)\cdot n(x) \geqslant 0$, then there exists $a'_0 \in A_0$ such that $-h(x,a'_0)\cdot n(x) >0$.

\medskip

\noindent \textbf{(HBB')} \quad For all $x \in \partial E$,  if $-h(x,a_0)\cdot n(x) \geqslant 0$ for all $a_0 \in A_0$, then $-h(x,a_0)\cdot n(x) > 0$ for all $a_0 \in A_0$.

\medskip

Let us now consider the  Hamilton-Jacobi-Bellman equation  associated to the optimal control problem, which turns out to be an elliptic fully nonlinear integro-differential equation on $[0,\,\infty)\times \bar E$
with nonlocal boundary conditions
\begin{align}
H^v(x,v(x),\nabla v(x)) &= 0  \qquad\quad \text{in} \,\, E,\label{HJB_E}\\
v(x) &= F^v(x)\quad \text{on} \,\, \partial E,\label{HJB_partialE}
\end{align}
where
\begin{align*}
H^\psi(z,u,p) &:= \sup_{a_0 \in A_0}\left\{ \delta\,u  - h(z,a_0) \cdot p -f(z,a_0)- \int_{E} (\psi(y)-\psi(z))\, \lambda(z, a_0)\,  Q(z,a_0, dy) \right\},\\
F^\psi(x)&:= \min_{a_\Gamma \in A_\Gamma} \left\{c(z,a_\Gamma) + \int_{E} \psi(y)\,  R(z,a_\Gamma, dy)\right\}.
\end{align*}
In the following the shorthand u.s.c. (resp. l.s.c.) stands for upper (resp. lower) semicontinuous.
\begin{Definition}\label{Def_HJB_viscosity_sol_2}
\begin{itemize}
\item[(i)] A bounded u.s.c. function $u$ on $\bar{E}$ is a viscosity subsolution of \eqref{HJB_E}-\eqref{HJB_partialE} if and only if, $\forall \phi \in \mathbb C^1_b(\bar{E})$, 
if $x_0\in \bar{E}$ is a global maximum of $u-\phi$ one has
\begin{eqnarray*}
&&H^u(x_0,u(x_0),\nabla\phi(x_0)) \leqslant 0 \,\,\, \qquad \qquad \qquad \qquad \qquad \text{if} \,\, x_0 \in E,\\
&&\min\{H^u(x_0,u(x_0),\nabla\phi(x_0)),u(x_0) - F^u(x_0)\}  \leqslant 0 \quad \text{if} \,\, x_0 \in \partial E.
\end{eqnarray*}
\item[(ii)] A bounded l.s.c. function $w$ on $\bar{E}$ is a viscosity supersolution of \eqref{HJB_E}-\eqref{HJB_partialE} if and only if, $\forall \phi \in \mathbb C^1_b(\bar{E})$, if $x_0\in \bar{E}$ is a global minimum of $w-\phi$ one has
\begin{eqnarray*}
&&H^w(x_0,w(x_0),\nabla\phi(x_0)) \geqslant 0 \,\,\, \qquad \qquad \qquad \qquad \qquad \text{if} \,\, x_0 \in E,\\
&&\max\{H^w(x_0,w(x_0),\nabla\phi(x_0)),w(x_0) - F^w(x_0)\}  \geqslant 0 \quad \text{if} \,\, x_0 \in \partial E.
\end{eqnarray*}
\item[(iii)] A viscosity solution of \eqref{HJB_E}-\eqref{HJB_partialE} is 
a continuous function which is both a viscosity subsolution and a viscosity supersolution of  \eqref{HJB_E}-\eqref{HJB_partialE}.
\end{itemize}
\end{Definition}
The following  theorem  collects the results of    Theorems   5.8 and 7.5 
 in \cite{Da-Fa}.
\begin{theorem}\label{Thm_HJB_unique_viscosity_sol}
Let \textup{\textbf{(H$\textup{h$\lambda$QR}$)}}, \textup{\textbf{(H$\textup{fc}$)}}, \textup{\textbf{(H$0$)}}, \textup{\textbf{(HBB)}} and \textup{\textbf{(HBB')}}  hold, and assume  that $A_0$, $A_\Gamma$ are compact. Let $V: E \rightarrow \R$ 
be the value function of the  PDMPs optimal control problem \eqref{Sec:PDP_value_function}.
 Then $V$ is a bounded and continuous function, and is	
the unique viscosity solution of \eqref{HJB_E}-\eqref{HJB_partialE}. 
\end{theorem}

\section{The randomized optimal control problem}
\label{Sec:PDP_Section_dual_control}

In the present section we formulate the randomized optimal control problem. First we introduce some notations. For every $a_0\in A_0$, we denote by $\phi(t, x,a_0)$ the unique solution to the ordinary differential equation
 $$
 \dot x(t) = h(x(t), a_0), \quad x(0)=x \in E.
 $$
 Notice that $\phi(t,x,a_0)$ coincides with the function $\phi^U(t,x)$, introduced in Section \ref{Sec:PDP_Sec_control_problem}, when $U(t)\equiv a_0$.
We also introduce two positive measures $\lambda_0$ and $\lambda_\Gamma$
 on $(A_0, \mathcal A_0)$ and $(A_\Gamma, \mathcal A_\Gamma)$, respectively, satisfying the following assumption:
 
 \medskip
 
 \noindent \noindent \textbf{(H\textup{$\lambda_0\lambda_\Gamma$})}\quad $\lambda_0$ and $\lambda_\Gamma$ are two finite positive measures on $(A_0, \mathcal A_0)$ and $(A_\Gamma, \mathcal A_\Gamma)$, respectively, with full topological support.
 
 \medskip
 
\noindent  For all $t \geq 0$, $(a_0,a_\Gamma) \in A_0 \times A_\Gamma$, let us  define
 \begin{eqnarray}
 \tilde \phi(t,x,a_0, a_\Gamma)&:=& (\phi(t,x,a_0), \, a_0, \, a_\Gamma),\quad x \in  \bar E, 
 \notag\\ 
 \tilde{\lambda}(x,a_0) &:=& \lambda(x,a_0)+ \lambda_0(A_0) + \lambda_\Gamma(A_\Gamma),\quad x \in  \bar E, \label{Sec:PDP_lambda_XI}\\ 
  \tilde{R}(x,a_0, a_\Gamma,dy\,db\,dc) &:=& R(x,a_\Gamma,dy)\,\delta_{a_0}(db)\,\delta_{a_\Gamma}(dc), \quad x \in \partial E,\notag 
 \end{eqnarray}
 and,  
 \begin{align*}
 &\tilde{Q}(x,a_0,a_\Gamma,dy\,db\,dc) :=\nonumber\\
 &\frac{\lambda(x,a_0)\,Q(x,a_0,dy)\,\delta_{a_0}(db)\,\delta_{a_\Gamma}(dc) + \lambda_0(db)\,\delta_{a_\Gamma}(dc)\,\delta_{x}(dy)+\lambda_\Gamma(dc)\,\delta_{a_0}(db)\,\delta_{x}(dy)}{\tilde{\lambda}(x,a_0)}, \quad x \in \bar E,\notag
 \end{align*}
 where, for any $F$ topological space, $\delta_a$ denotes the Dirac measure concentrated at some point  $a \in F$.

\subsection{State process}
 
Our purpose is  to construct a PDMP $(X,I, J)$  with enlarged state space $E \times A_0 \times A_\Gamma$ and  local characteristics $(\tilde \phi,\tilde \lambda, \tilde Q, \tilde R)$. This can be done in a canonical way, proceeding as in Section \ref{Sec:PDP_Sec_control_problem}. By an abuse of notation, we use the same symbols as in Section \ref{Sec:PDP_Sec_control_problem}. So, in particular, we  define  $\Omega'$ as the set of sequences $\omega'=(t_n, e_n, a_n^0,  a_n^{\Gamma})_{n \geq 1}$ contained in $((0,\,\infty) \times E \times A_0 \times A_\Gamma) \cup \{ (\infty, \Delta, \Delta', \Delta'')\}$, where $\Delta \notin E$, $\Delta' \notin A_0$, $\Delta'' \notin A_\Gamma$  are isolated points respectively adjoined to $E$, $A_0$ and $A_\Gamma$. In the sample space $\Omega = \Omega' \times E \times A_0 \times A_{\Gamma}$ we define the random variables $T_0(\omega)=0$, $E_0(\omega)=x$,  $A^0(\omega)=a_0$,  $A^\Gamma(\omega)=a_\Gamma$, and the sequence of random variables   $T_n : \Omega \rightarrow (0,\,\infty]$,
 $E_n : \Omega \rightarrow E \cup \{\Delta\}$, 
 $A^0_n : \Omega \rightarrow A_0 \cup \{\Delta'\}$,  
 $A^{\Gamma}_n : \Omega \rightarrow A_{\Gamma} \cup \{\Delta''\}$, for $n \geq 1$, by setting $T_n (\omega )= t_n$, 
$E_n(\omega)=e_n$,
$A^0_n(\omega)=a^0_n$,
$A^\Gamma_n(\omega)=a^\Gamma_n$, with
$T_\infty  (\omega )= \lim_{n\to\infty} t_n$.
 Then, we  define the process $(X,I,J)$ on $(E \times A_0 \times A_{\Gamma}) \cup \{\Delta, \Delta', \Delta''\}$ as 
 \begin{align*}
 (X,I,J)_t&=
 \left\{
 \begin{array}{ll}
 (\phi(t-T_n, E_n,A_n^0),A_n^0, A_n^\Gamma)& \quad \textup{if}\,\, T_n \leq t < T_{n+1},\,\,\textup{for}\,\,n \in \N,\\
 (\Delta,\Delta', \Delta'')&  \quad \textup{if}\,\,t \geq T_{\infty}.
 \end{array}	
 \right.
 \end{align*}
We also define the  random measure $p$  on $(0,\,\infty) \times E \times A$ as
 \begin{equation}\label{Sec:PDP_p_dual}
 p(ds\,dy\,db\,dc)= \sum_{n\in \N} \one_{\{T_n,E_n, A^0_n, A_n^{\Gamma}\}}(ds\,dy\,db\,dc),
 \end{equation}
and, for all $t\geq 0$, we  introduce   the $\sigma$-algebra
 $\mathcal G_t=\sigma(p((0,s] \times G)\,:\, s\in (0,t], G\in\mathcal E \otimes \mathcal A_0 \times \mathcal A_{\Gamma})$, 
and the
 $\sigma$-algebra $\mathcal F_t$ generated by $\mathcal F_0$ and $\mathcal G_t$, where
 $\mathcal F_0=\mathcal E\otimes \mathcal A_0 \otimes \mathcal A_{\Gamma} \otimes \{\emptyset,\Omega'\}$.
 We still denote by $\F = (\mathcal F_t)_{t \geq 0}$ and $\mathcal P$ the corresponding filtration and predictable $\sigma$-algebra.

 Given any starting point $(x,a_0, a_{\Gamma}) \in E \times A^0\times A^{\Gamma}$,
 by  Proposition 2.1 in \cite{BandiniPDMPsNoBordo}, there exists a unique
 probability measure    on $(\Omega, \mathcal F_{\infty})$, denoted by $\P^{x,a_0,a_\Gamma}$, such that its restriction to  $\mathcal F_0$ is $\delta_{(x,a_0,a_\Gamma)}$  and the $\mathbb F$-compensator  of the measure $p(ds\,dy\,db\,dc)$ under $\P^{x,a_0,a_\Gamma}$ is the random measure
 $$
\tilde p(ds\,dy\,db\,dc)= \sum_{n \in \N}\one_{[T_n,\,T_{n+1})}(s)\,\Lambda( \phi(s-T_n, E_n,A_n^0),A_n^0,A_n^\Gamma, dy\,db\,dc)\,d A_s,
 $$
 where, for all $(x,a_0,a_\Gamma)\in E \times A_0 \times A_\Gamma$, 
 \begin{align*}
 &\Lambda(x,a_0,a_\Gamma,dy\,db\,dc)= 
\tilde{Q}(x,a_0,a_\Gamma,dy\,db\,dc) \,\one_{x \in E}+
 \tilde{R}(x,a_0,a_\Gamma,dy\,db\,dc)\,\one_{x \in \Gamma},
 \end{align*}
 and $A_s$ is the increasing, predictable process such that, for any $s \geq 0$, 
 \begin{align*}
d A_s(\omega) = \tilde \lambda(X_{s-}(\omega),I_{s-}(\omega))\,\one_{X_{s-}(\omega) \in E}\,ds + \one_{X_{s-}(\omega) \in \Gamma}\,d p^{\ast}_s(\omega).
 \end{align*}
 In particular,    
\begin{align}
&\Delta A_t(\omega) = \tilde p (\omega, \{t\} \times E \times A_0 \times A_\Gamma) = \one_{X_{t-}(\omega) \in \partial E}\label{DeltaA},\\
 &\tilde p (\omega, \{t\} \times dy\, db\,dc) = \tilde R(X_{t-}(\omega), I_{t-}(\omega), J_{t-}(\omega),dy \,db\,dc)\,\Delta A_t(\omega).\label{tildep_t}
\end{align}
 \begin{Remark}
The  $\mathbb F$-compensator  of the measure $p(ds\,dy\,db\,dc)$ under $\P^{x,a_0,a_\Gamma}$ can  be decomposed as 
$\tilde p(\omega, ds\,dy\,db\,dc) = \phi_{\omega,t}(dy\,db\,dc)\, d A_s(\omega)$,
where
\begin{equation}\label{phi_omegat}
\phi_{\omega,t}(dy\,db\,dc):= \Lambda(X_{t-}(\omega),I_{t-}(\omega),J_{t-}(\omega),dy\,db\,dc).
\end{equation}
 \end{Remark}
The process $(X,I, J)$ is Markovian on $[0,\,\infty)$ with respect to $\F$.
 For
 every 
 real-valued functions $\varphi$  defined on $E \times A_0 \times A_\Gamma$, 
 we define
 \begin{align*}
 &\mathcal{L}\varphi(x,a_0,a_\Gamma) :=h(x,a_0)\cdot \nabla_x \varphi(x,a_0,a_\Gamma) + \int_{E}(\varphi(y,a_0,a_\Gamma) - \varphi(x,a_0,a_\Gamma))\, \lambda(x,a_0)\,Q(x,a_0,dy)\\
 &+\int_{A_0}(\varphi(x,b,a_\Gamma) - \varphi(x,a_0,a_\Gamma))\, \lambda_0(db)+\int_{A_\Gamma}(\varphi(x,a_0,c) - \varphi(x,a_0,a_\Gamma))\, \lambda_\Gamma(dc),\quad x \in E,\\
&\mathcal{G}\varphi(x,a_0,a_\Gamma) := \int_{E}(\varphi(y,a_0,a_\Gamma) - \varphi(x,a_0,a_\Gamma))\,R(x,a_\Gamma,dy),\quad x \in \partial E.
 \end{align*}
 From Theorem 26.14 in \cite{Da} it follows that  $\mathcal{L}$ is the extended generator of the process $(X,I,J)$ and $\mathcal{G}\varphi=0$ if and only if $\varphi$ belongs to the domain of $\mathcal{L}$.

 \subsection{The control problem}\label{Sec:PDP_Section_dual_optimal_control}

The class of admissible control maps  is the set $\Vc= \mathcal V_0 \otimes \mathcal V_\Gamma$, where 
 \begin{align*}
 \mathcal{V}_0 = \{ \nu^0: \Omega \times [0,\,\infty) \times A _0\rightarrow (0,\,\infty)\,\, \mathcal{P}\otimes \mathcal A_0\text{-measurable and bounded}\},\\
  \mathcal{V}_\Gamma = \{ \nu^\Gamma: \Omega \times [0,\,\infty) \times A_\Gamma \rightarrow (0,\,\infty)\,\, \mathcal{P}\otimes \mathcal A_\Gamma\text{-measurable and bounded}\}.
 \end{align*}
 For every $\bm\nu=(\nu^0,\nu^\Gamma) \in \mathcal{V}$, 
 we define
 \begin{align*}
  &\tilde{\lambda}^{\bm\nu}(t,x,a_0) := \lambda(x,a_0)+ \int_{A_0}\nu^0_t(b)\,\lambda_0(db) + \int_{A_\Gamma}\nu^\Gamma_t(c)\,\lambda_\Gamma(dc),\\
 &\tilde{Q}^{\bm\nu}(t,x,a_0,a_\Gamma,dy\,db\,dc) :=\notag\\
 &\frac{\lambda(x,a_0)\,Q(x,a_0,dy)\,\delta_{a_0}(db)\,\delta_{a_\Gamma}(dc) + \nu^0_t(b)\,\lambda_0(db)\,\delta_{a_\Gamma}(dc)\,\delta_{x}(dy)+\nu^\Gamma_t(c)\,\lambda_\Gamma(dc)\,\delta_{a_0}(db)\,\delta_{x}(dy)}{\tilde{\lambda}^{\bm\nu}(t,x,a_0)},\notag \\
&\Lambda^{\bm\nu}(t,x,a_0,a_\Gamma,dy\,db\,dc)=
\tilde{Q}^{\bm\nu}(t,x,a_0,a_\Gamma,dy\,db\,dc) \,\one_{x \in E}+
 \tilde{R}(x,a_0,a_\Gamma,dy\,db\,dc)\,\one_{x \in \partial E}.\notag
 \end{align*}
Then,  for every $\bm\nu \in \mathcal{V}$, we consider the predictable random measure
 \begin{equation}\label{Sec:PDP_dual_comp}
\tilde p^{\bm\nu}(ds\,dy\,db\,dc)= \sum_{n \in \N}\one_{[T_n,\,T_{n+1})}(s)\,\Lambda^{\bm\nu}(s,\phi(s-T_n, E_n,A_n^0),A_n^0,A_n^\Gamma, dy\,db\,dc)\,d A^{\bm\nu}_s,
\end{equation}
 where $A^{\bm\nu}$ is the increasing and predictable process given by 
  \begin{align*}
dA_s^{\bm\nu} = \tilde\lambda^{\bm\nu}(s,X_{s-},I_{s-})\,\one_{X_{s-} \in E}\,ds + \one_{X_{s-} \in \partial E}\,d p^\ast_s.
  \end{align*}
In what follows we will denote $q=p-\tilde{p}$ and $q^{\bm\nu}= p-\tilde{p}^{\bm\nu}$.
We have the following important result.

\begin{proposition}\label{P_Prob_infty}
 	Let assumptions {\bf (Hh$\lambda$QR)}
 	 and  {\bf (H$\lambda_0\lambda_\Gamma$)} hold. Then, for every  $(x,a_0,a_\Gamma)\in E\times A_0 \times A_\Gamma$  and
 	$\bm\nu \in \mathcal{V}$, there exists a unique probability
 	$\P^{x,a_0,a_\Gamma}_{\bm\nu}$ on $(\Omega, \mathcal F_{\infty})$, under which the random measure  $\tilde{p}^{\bm\nu}$ in \eqref{Sec:PDP_dual_comp} is the compensator of the measure $p$ in \eqref{Sec:PDP_p_dual} on $(0,\,\infty)\times E \times A_0 \times A_\Gamma$.
\end{proposition}
\proof
The proof of Proposition \ref{P_Prob_infty} is postponed in Appendix \ref{S:AppP}, see Proposition \ref{Sec:PDP_lemma_P_nu_martingale}.
\endproof

For every $(x,a_0,a_\Gamma) \in E \times A_0 \times A_\Gamma$, the randomized optimal control problem consists in minimizing over all $\bm\nu \in \mathcal V$ the cost functional (we denote by $\mathbb{E}_{\bm\nu}^{x,a_0,a_\Gamma}$ the expectation operator under $\P_{\bm\nu}^{x,a_0,a_\Gamma}$)
$$
 J(x,a_0,a_\Gamma,\bm\nu) := \spernuxa{ \int_{(0, \infty)} e^{-\delta\,s}\, f(X_{s},I_s)\,ds+ \int_{(0,\infty)} e^{-\delta\,s}\, c(X_{s-},J_{s-})\,d p^\ast_s},
$$
The value function is given by
\begin{equation}\label{Sec:PDP_dual_value_function}
 V^{\ast}(x,a_0,a_\Gamma) := \inf_{ \bm\nu \in \mathcal{V}} J(x,a_0,a_\Gamma,\bm\nu).
\end{equation}

\section{Constrained BSDEs and probabilistic representation of $V^\ast$}
\label{Sec:PDP_Sec_ConstrainedBSDE}

In the present section we introduce a BSDE with two sign constraints on its martingale part, that will provide a probabilistic representation formula for the value function $V^*$ in \eqref{Sec:PDP_dual_value_function}. The main novelty with respect to previous literature is that our BSDE is driven by a non quasi-left-continuous random measure. For such an equation, the proof of existence and uniqueness is a difficult task, and 
counterexamples  can be obtained even in simple cases, see \cite{CFJ}.  Only recently,  some results in the unconstrained case have been obtained in this context, see \cite{CohenElliottPearce}, \cite{CohenElliott}, \cite{BandiniBSDE}. In order to have an existence and uniqueness result for our BSDE, we have to impose the following additional assumption on $p^\ast$.

\medskip

\noindent \textbf{(H$0'$)} \quad For any $(x,a_0,a_\Gamma) \in E \times A_0 \times A_\Gamma$, $t  \in \R_+$, $\beta >0$, there exists some  $\bar C_\beta(t) < \infty$, only depending  on $t$ and $\beta$, such that 
$\E^{x, a_0,a_\Gamma}\left[(1 + p^\ast_t)\,(1+ \beta)^{p^\ast_t}\right]\leqslant \bar C_\beta(t)
$.

\medskip
 
Now, we introduce some notations. Firstly, for any $\beta\geq0$, given a predictable increasing process $A$, we  denote by $\mathcal E^\beta$ the Dol\'eans-Dade exponential of the process $\beta A$, given by
\begin{equation*}
\mathcal{E}_t^{\beta} \ = \ e^{\beta\,A_t}\prod_{0 < s \leq t}(1 + \beta\,\Delta A_s)\,e^{-\beta\,\Delta\,A_s}.
\end{equation*}
In particular, $d \mathcal{E}_{s}^{\beta} = \beta\,\mathcal{E}_{s-}^{\beta}\,dA_s$, $\mathcal{E}_{s}^{\beta} \geq 1$.
 \begin{Remark}
Given a  c\`adl\`ag process $C$, It\^o's formula applied to $\mathcal{E}_{s}^{\beta}\,|C_s|^2$ yields 
\begin{align*}
d(\mathcal{E}_{s}^{\beta}\,|C_s|^2) 
&= 2\,\mathcal{E}_s^{\beta}\,C_{s-}\,d C_s + \mathcal{E}_s^{\beta}\,(\Delta C_s)^2+\beta \, \mathcal{E}_{s-}^{\beta}\,|C_{s-}|^2\,d A_s \notag \\
&= 2\,\mathcal{E}_s^{\beta}\,C_{s-}\,d C_s + \mathcal{E}_s^{\beta}\,(\Delta C_s)^2+\beta \, \mathcal E_s^\beta\,(1+\beta\Delta A_s)^{-1}\,|C_{s-}|^2\,d A_s,
\end{align*}
where the last equality follows from the fact that $\mathcal{E}_{s-}^{\beta} = \mathcal{E}_{s}^{\beta}\,(1+ \beta\,\Delta A_s)^{-1}$.
 \end{Remark}
For any $(x,a_0,a_\Gamma) \in E \times A_0 \times A_\Gamma$, and  	$\beta\geq0$, we  introduce the following sets.
 \begin{itemize}
 	\item $\textbf{L}^{\textbf{2}}_{\textbf{x},\textbf{a}_0, \textbf{a}_\Gamma}(\mathcal{F}_\tau)$,  the set of $\mathcal{F}_\tau$-measurable random variables $\xi$ such that $\sperxa{
|\xi|^2} < \infty$; here  $\tau \geqslant 0$ is an $\mathbb F$-stopping time.
 	\item $\textbf{S}^\infty$  the set of real-valued c\`adl\`ag 	adapted processes $Y = (Y_t)_{t \geqslant 0}$ which are uniformly bounded.
 	\item $\textbf{L}_{\textbf{x},\textbf{a}_0, \textbf{a}_\Gamma}^{\textbf{2},\beta}(\textup{p}^\ast;\textbf{0,\,T})$, $T >0$,  the set of real-valued progressive processes $Y= (Y_t)_{0 \leqslant t \leqslant T}$ such that
 	\begin{displaymath}
 	||Y||^2_{\textbf{L}^{\textbf{2},\beta}_{\textbf{x},\textbf{a}_0, \textbf{a}_\Gamma}(\textup{p}^\ast;\textbf{0,\,T})}:=\sperxa{\int_{(0,\,T]} \mathcal E_t^\beta\, |Y_{t-}|^2\,d A_t} < \infty.
 	\end{displaymath}
 	\item 
 	$\mathcal{G}_{\textbf{x},\textbf{a}_0, \textbf{a}_\Gamma}^{\textbf{2},\beta}(\textup{q};\textbf{0,\,T})$, $T >0$,  the set of $\mathcal{P}_T\otimes \mathcal E \otimes \mathcal A_0 \otimes \mathcal A_\Gamma$-measurable maps $Z: \Omega \times [0,\,T] \times E \times A_0 \times A_\Gamma \rightarrow \R$ such that
 	\begin{align*}
 	||Z||^2_{\mathcal{G}_{\textbf{x},\textbf{a}_0, \textbf{a}_\Gamma}^{\textbf{2}, \beta}(\textup{q};\textbf{0,\,T})}&:=
 	\ \sperxa{\int_{(0,\,T]}\mathcal E_t^\beta\int_{E\times A_0 \times A_\Gamma} \big|Z_t(y,b,c) - \hat Z_t \,\one_{K}(t)\big|^2\,\tilde p(dt\,dy\,db\,dc) 
	}
		\end{align*}
		is finite, where
\begin{equation}\label{Zhat}
		\hat Z_t \ := \ \int_{E\times A_0 \times A_\Gamma}  Z_t(y,b,c)\,\tilde{p}(\{t\} \times dy\,db\,dc), \qquad 0\leq t\leq T.
		\end{equation}
We also define $\mathcal{G}^{\textbf{2},\beta}_{{\textbf{x},\textbf{a}_0, \textbf{a}_\Gamma},\textbf{loc}}(\textup{q}) :=\cap_{T >0} \,\mathcal{G}^{\textbf{2}, \beta}_{\textbf{x},\textbf{a}_0, \textbf{a}_\Gamma}(\textup{q};\textbf{0,\,T})$. When $\beta=0$, we simply write $\mathcal{G}_{\textbf{x},\textbf{a}_0, \textbf{a}_\Gamma}^{\textbf{2}}(\textup{q};\textbf{0,\,T})$ and $\mathcal{G}^{\textbf{2}}_{{\textbf{x},\textbf{a}_0, \textbf{a}_\Gamma},\textbf{loc}}(\textup{q})$ in place of $\mathcal{G}_{\textbf{x},\textbf{a}_0, \textbf{a}_\Gamma}^{\textbf{2},\beta}(\textup{q};\textbf{0,\,T})$ and $\mathcal{G}^{\textbf{2},\beta}_{{\textbf{x},\textbf{a}_0, \textbf{a}_\Gamma},\textbf{loc}}(\textup{q})$, respectively (for equivalent notions of the $\mathcal{G}_{\textbf{x},\textbf{a}_0, \textbf{a}_\Gamma}^{\textbf{2}}(\textup{q};\textbf{0,\,T})$ norm see Lemma \ref{L_simplifiedNorm}).
 	\item $\textbf{L}^{\textbf{2}}(\lambda_0)$ (resp.  $\textbf{L}^{\textbf{2}}(\lambda_\Gamma)$),
 	the set of $\mathcal{A}_0$-measurable maps $\psi: A_0 \rightarrow \R$ (resp. $\mathcal{A}_{\Gamma}$-measurable maps $\psi: A_\Gamma \rightarrow \R$) such that
 	\begin{eqnarray*}
 		|\psi|^2_{\textbf{L}^{\textbf{2}}(\lambda_0)}:=\int_{A_0} |\psi(b)|^2 \, \lambda_0(db)< \infty\quad \left(\textup{resp.}\quad |\psi|^2_{\textbf{L}^{\textbf{2}}(\lambda_\Gamma)}:=\int_{A_\Gamma} |\psi(c)|^2 \, \lambda_\Gamma(dc)< \infty\right).
 	\end{eqnarray*}	
 	 \item  $\textbf{L}^{\textbf{2}}(\phi_{\omega,t})=\textbf{L}^2(E \times A_0 \times A_\Gamma, \mathcal E \otimes \mathcal A_0\otimes \mathcal A_\Gamma,\phi_{\omega,t}(dy\,db\,dc))$,
for any $(\omega,t) \in \Omega \times \R_+$, the set of $\mathcal E \otimes \mathcal A_0\otimes \mathcal A_\Gamma$-measurable maps $\zeta: E\times A_0 \times A_{\Gamma} \rightarrow \R$ such that  
 	 \begin{eqnarray*}
 		|\zeta|^2_{\textbf{L}^{\textbf{2}}(\phi_{\omega,t})}:=\int_{E \times A_0 \times A_\Gamma} |\zeta(y,b,c)|^2 \, \phi_{\omega,t}(dy\,db\,dc)< \infty,
 	\end{eqnarray*}
 	where $\phi_{\omega,t}(dy\,db\,dc)$ is the random measure introduced in \eqref{phi_omegat}. 
 	\item $\textbf{K}^{\textbf{2}}_{\textbf{x},\textbf{a}_0, \textbf{a}_\Gamma}(\textbf{0,\,T})$, $T>0$,  the set of nondecreasing c\`adl\`ag predictable  processes $K = (K_t)_{0 \leqslant t \leqslant T}$ such that $K_0 = 0$ and $\sperxa{
 	|K_T|^2}< \infty$. We also define $\textbf{K}^{\textbf{2}}_{{\textbf{x},\textbf{a}_0, \textbf{a}_\Gamma},\textbf{loc}} :=\cap_{T >0} \textbf{K}^{\textbf{2}}_{\textbf{x},\textbf{a}_0, \textbf{a}_\Gamma}(\textbf{0,\,T})$.
 \end{itemize}
  We aim at  studying the following family of BSDEs with partially nonnegative jumps over an infinite horizon, parametrized by $(x,a_0,a_\Gamma)$: 
 ${\mathbb{P}}^{x,a_0,a_\Gamma}$-a.s.,
 \begin{align}
 & Y^{x,a_0,a_{\Gamma}}_{s} = Y^{x,{a_0,a_{\Gamma}}}_T - \delta \,\int_{(s,\,T]}Y^{x,{a_0,a_{\Gamma}}}_r\,dr + \int_{(s,\,T]}f(X_r,I_r)\,dr + \int_{(s,\,T]}c(X_{r-},J_{r-})\,d p^\ast_r \nonumber\\
 &- \int_{(s,\,T]}\int_{A_0}Z^{x,a_0,a_{\Gamma}}_r(X_r,\,b,\,J_r)\, \lambda_0(db)\,dr 
 - \int_{(s,\,T]}\int_{A_{\Gamma}}Z^{x,a_0,a_{\Gamma}}_r(X_r,\,I_r,\,c)\, \lambda_{\Gamma}(dc)\,dr
 \nonumber\\
 &-\big(K^{x,a_0,a_\Gamma}_T - K^{x,a_0,a_\Gamma}_s\big)- \int_{(s,\,T]}\int_{E \times  A_0 \times A_{\Gamma}}Z^{x,a_0,a_{\Gamma}}_r(y,\,b,\,c)\, q(dr\,dy\, db\, dc),\quad  0\leqslant s \leqslant T<\infty,\label{Sec:PDP_BSDE}
 \end{align}
 with the constraints
 \begin{align}
 Z^{x,a_0,a_\Gamma}_s(X_{s-},b, J_{s-})\geqslant 0,\quad d\P^{x,a_0,a_\Gamma} &\lambda_0(db)\textup{ -a.e. on } [0,\,\infty)\times \Omega \times A_0,\label{Sec:PDP_BSDE_constraint1}\\
  Z^{x,a_0,a_\Gamma}_s(X_{s-},I_{s-},c)\geqslant 0\quad d\P^{x,a_0,a_\Gamma} &\lambda_{\Gamma}(dc)\textup{ -a.e. on } [0,\,\infty)\times \Omega \times A_\Gamma.\label{Sec:PDP_BSDE_constraint2}
 \end{align}
We look for a \emph{maximal solution} $(Y^{x,a_0,a_\Gamma},Z^{x,a_0,a_\Gamma},K^{x,a_0,a_\Gamma})\in \textbf{S}^{\infty}\times \mathcal{G}^{\textbf{2}}_{\textbf{x},\textbf{a}_0,\textbf{a}_\Gamma,\textbf{\text{loc}}}(\textup{q})\times \textbf{K}^{\textbf{2}}_{\textbf{x},\textbf{a}_0,\textbf{a}_\Gamma,\textbf{\text{loc}}}$ to \eqref{Sec:PDP_BSDE}-\eqref{Sec:PDP_BSDE_constraint1}-\eqref{Sec:PDP_BSDE_constraint2},
 in the sense that for any other solution $(\tilde{Y}, \tilde{Z},\tilde{K})\in \textbf{S}^{\infty}\times \mathcal{G}^{\textbf{2}}_{\textbf{x},\textbf{a}_0,\textbf{a}_\Gamma,\textbf{\text{loc}}}(\textup{q})\times \textbf{K}^{\textbf{2}}_{\textbf{x},\textbf{a}_0,\textbf{a}_\Gamma,\textbf{\text{loc}}}$ to \eqref{Sec:PDP_BSDE}-\eqref{Sec:PDP_BSDE_constraint1}-\eqref{Sec:PDP_BSDE_constraint2}, we  have $Y_t^{x,a_0,a_\Gamma} \geqslant \tilde{Y}_t$, $\P^{x,a_0,a_\Gamma}$-a.s., for all $t\geqslant 0$.

We can now state the main result of this section.
    
\begin{theorem}\label{Sec:PDP_Thm_ex_uniq_maximal_BSDE}
 	Under assumptions  \textup{\textbf{(H$\textup{h$\lambda$QR}$)}}, {\bf (H$0$)}, {\bf (H$0'$)},  \textup{\textbf{(H$\lambda_0\lambda_\Gamma$)}} and  \textup{\textbf{(H$\textup{fc}$)}},
 	for every  $(x,a_0,a_\Gamma) \in E \times A_0 \times A_\Gamma$, 
 	there exists a unique maximal solution $(Y^{x,a_0,a_\Gamma},Z^{x,a_0,a_\Gamma},K^{x,a_0,a_\Gamma})\in \textup{\textbf{S}}^{\infty}\times \mathcal{G}^{\textbf{2}}_{\textbf{x},\textbf{a}_0,\textbf{a}_\Gamma,\textbf{\textup{loc}}}(\textup{q})\times \textup{\textbf{K}}^{\textbf{2}}_{\textbf{x},\textbf{a}_0,\textbf{a}_\Gamma,\textbf{\textup{loc}}}$  to the constrained BSDE \eqref{Sec:PDP_BSDE}-\eqref{Sec:PDP_BSDE_constraint1}-\eqref{Sec:PDP_BSDE_constraint2}. Moreover, $Y^{x,a_0,a_\Gamma}$  has the explicit representation:
 	\begin{equation}\label{Sec:PDP_rep_Y}
 	Y_s^{x,a_0,a_\Gamma} = \essinf_{\bm\nu \in \mathcal{V}}\spernuxa{\int_{(s,\infty)} e^{-\delta (r-s)}\,f(X_{r},I_{r})\, dr +\int_{(s,\infty)} e^{-\delta (r-s)}\,c(X_{r-},J_{r-})\, d p^\ast_r\Big| \mathcal{F}_s},
 	\end{equation}
for all $s\geqslant0$. In particular, setting $s=0$, 
 	we have the following representation formula for the value function of the randomized control problem:
 	\begin{equation}\label{Sec:PDP_Vstar_Y0}
 	V^{\ast}(x,a_0,a_\Gamma)
 	= Y_0^{x,a_0,a_\Gamma}, \quad (x,a_0,a_\Gamma)\in E \times A_0 \times A_\Gamma.
 	\end{equation}
 \end{theorem}
 \proof
We start by considering, for every $T>0$,  the family of penalized BSDEs on $[0,\,T]$  with zero  terminal condition at time $T$, parametrized by the integer $n \geqslant 1$: $\P^{x,a_0,a_\Gamma}$-a.s.
 \begin{align}\label{Sec:PDP_BSDE_penalized_T}
 Y_{s}^{T,n,x,a_0,a_\Gamma}&=  -\delta\int_{(s,T]} Y_r^{T,n,x,a_0,a_\Gamma}\,dr + \int_{(s,\,T]} f(X_{r},I_{r})\, dr+ \int_{(s,\,T]} c(X_{r-},J_{r-})\, d p^\ast_r\nonumber\\
 & -n\int_{(s,T]}\int_{A_0} [Z_r^{T,n,x,a_0,a_\Gamma}(X_r,b,J_r)]^-\,\lambda_0(db)\,dr - \int_{(s,\,T]} \int_{A_0} Z_{r}^{T,n,x,a_0,a_\Gamma}(X_{r},b,J_r) \,\lambda_{0}(db)\, dr \nonumber\\
  & -n\int_{(s,T]}\int_{A_\Gamma} [Z_r^{T,n,x,a_0,a_\Gamma}(X_r,I_r,c)]^-\,\lambda_\Gamma(dc)\,dr - \int_{(s,\,T]} \int_{A_{\Gamma}} Z_{r}^{T,n,x,a_0,a_\Gamma}(X_{r},I_r,c) \,\lambda_{\Gamma}(dc)\, dr \nonumber\\
 &- \int_{(s,\,T]} \int_{E \times A_0 \times A_{\Gamma}} Z_{r}^{T,n,x,a_0,a_\Gamma}(y,b,c) \, q(dr\,dy\,db\,dc),\quad 0\leqslant s\leqslant T,
 \end{align}
where $[z]^- =\max(-z,0)$ is  the negative part of $z$. Our aim is to exploit equation \eqref{Sec:PDP_BSDE_penalized_T} in order to construct the maximal solution $(Y^{x,a_0,a_\Gamma},Z^{x,a_0,a_\Gamma},K^{x,a_0,a_\Gamma})$, studying the limit of $(Y^{T,n}, Z^{T,n}) = (Y^{T,n,x,a_0,a_\Gamma},Z^{T,n,x,a_0,a_\Gamma})$ firstly as $T\rightarrow\infty$, and then as $n\rightarrow\infty$. Before analyzing the asymptotic behavior of $(Y^{T,n}, Z^{T,n})$, we need to prove the existence of a unique solution to equation \eqref{Sec:PDP_BSDE_penalized_T}. This is indeed a consequence of Theorem 4.1 in \cite{BandiniBSDE}. As a matter of fact, notice that equation \eqref{Sec:PDP_BSDE_penalized_T} can be rewritten as: $\P^{x,a_0,a_\Gamma}$-a.s.,
 	\begin{align}\label{Sec:PDP_BSDE_penalized_TBIS_comp}
 	 Y_{s}^{T,n,x,a_0,a_\Gamma}
 	&=\int_{(s,\,T]} \tilde{f}^n(r-s,X_{r-},I_{r-},\,J_{r-},\, Z_r^{T,n,x,a_0,a_\Gamma})\, d A_r\nonumber\\ 
 	&- \int_{(s,\,T]} e^{-\delta (r-s)}\int_{E \times A_0\times A_{\Gamma}}  Z_{r}^{T,n,x,a_0,a_\Gamma}(y,b,c) \, q(dr\,dy\,db\,dc),\quad s \in [0,\,T],
 	\end{align}
where
\begin{equation}\label{Sec:PDP_tildef^n}
\tilde {f}^n(t,x,a_0, a_{\Gamma},\zeta)
 	:= e^{-\delta\,t}\,f^n(x,a_0,\zeta(x,\cdot,a_\Gamma),\zeta(x,a_0,\cdot))\,\one_{x \in E} +e^{-\delta\,t}\,c(x,a_{\Gamma})\,\one_{x \in \partial E},
\end{equation}
with
$f^n(x,a_0,\psi,\varphi) 
 	:= f(x,a_0)  
 	 -\int_{A_0}\left\{ n\,[\psi(b)]^{-} + \psi(b)\right\} 
 	\lambda_{0}(db)
 	-\int_{A_\Gamma}\left\{ n\,[\varphi(c)]^{-} + \varphi(c)\right\}
 	\lambda_{\Gamma}(dc)$.
Under  assumptions \textup{\textbf{(H$\lambda_0\lambda_\Gamma$)}} and \textup{\textbf{(H$\textup{fc}$)}}, there exists a constant $L_n$, depending only on $n$, such that
 	\begin{align}\label{fn_Lipschitzianity}
 	|f^n(x,a_0,\psi,\varphi)-f^n(x,a_0,\psi',\varphi')|\leqslant L_n \left(|\psi-\psi'|_{\textup{\textbf{L}}^{\textbf{2}}(\lambda_0)} + |\varphi-\varphi'|_{\textup{\textbf{L}}^{\textbf{2}}(\lambda_\Gamma)} \right),
 	\end{align}	
for every $(x,a_0)\in  E \times A_0 $, $\psi,\,\psi' \in \textup{\textbf{L}}^{\textbf{2}}(\lambda_0)$, $\varphi,\,\varphi' \in \textup{\textbf{L}}^{\textbf{2}}(\lambda_\Gamma)$. 
Then, one can easily show that 
\begin{align*}
		&|\tilde{f}^n(t,X_{t-}(\omega),I_{t-}(\omega),J_{t-}(\omega), \zeta') - \tilde{f}^n(t,X_{t-}(\omega),I_{t-}(\omega),J_{t-}(\omega), \zeta)| \leq 
		\notag \\
		& 2 \, L_n \bigg(\int_{E\times A_0 \times A_\Gamma} \bigg|\tilde \zeta(y,b,c)- \Delta A_t(\omega)\int_{E\times A_0 \times A_\Gamma}\,\tilde\zeta(\bar y,\bar b,\bar c)\,\phi_{\omega,t}(d\bar y\,d\bar b\,d\bar c)\bigg|^2\phi_{\omega,t}(dy\,db\,dc)\bigg)^{1/2}, 
		\end{align*}
		for all $(\omega,t)\in\Omega\times[0,T]$,
		 $\zeta,\zeta'\in \textbf{L}^{\textbf{2}}(\phi_{\omega,t})
		 $, with	 $\tilde \zeta = \zeta - \zeta'$.
Finally, setting  $c_1(T) = (M_f^2 \vee M_c^2)||\tilde \lambda||_\infty \,T \,e^{\beta ||\tilde \lambda||_\infty T}$, $c_2 = (M_f^2 \vee M_c^2)$, we have 
	\begin{align*}
	&\sperxa{(1 + \sum_{0<t\leq T} |\Delta A_t|^2)\int_{(0,\,T]} \mathcal{E}_t^{\beta} \left|\tilde{f}^n(t,X_{t-},I_{t-},J_{t-},0)\right|^2 d A_t}\\
	&\leq  c_1(T) (1 + \sperxa{p^\ast_T})+c_2\,	\sperxa{(1 + p^\ast_T)\,  (1+ \beta)^{p^\ast_T}},
	\end{align*}
which is finite by 
\eqref{conv_exp_sper_pstar} and hypothesis \textbf{(H$0'$)}.
We are therefore  in condition to apply   Theorem 4.1 in \cite{BandiniBSDE}.
Setting   
$$
\beta_0^n:=\frac{2\, (L_n + \varepsilon)^2}{1-\varepsilon}, \quad \varepsilon \in (0,1),
$$
  we deduce that there exists   of a unique solution $(Y^{T,n,x,a_0,a_\Gamma},Z^{T,n,x,a_0,a_\Gamma})\in
	\textup{\textbf{L}}_{\textbf{x},\textbf{a}_0, \textbf{a}_\Gamma}^{\textbf{2},\beta}(\textup{p}^\ast;\textup{\textbf{0,\,T}})
	\times \mathcal{G}^{\textbf{2},\beta}_{{\textbf{x},\textbf{a}_0, \textbf{a}_\Gamma}}(\textup{q};\textup{\textbf{0,\,T}})$ to equation \eqref{Sec:PDP_BSDE_penalized_T} 
	for  $\beta \geq \beta_0^n$. Notice that the Lischitz constant of $\tilde f^n$ with respect to $Y$, that we will denote  $L_y$, is identically zero. So, in particular, the technical assumption of Theorem 4.1 in \cite{BandiniBSDE}, that is the existence  of  $\varepsilon\in(0,1)$ such that (in our framework, $\Delta A_t = \one_{X_{t-} \in \partial E}$)
\begin{equation*}
2\,L_y^2\,|\Delta A_t|^2 \ \leq \ 1 - \varepsilon, \qquad \P\text{-a.s.},\,\forall\,t\in[0,T],
\end{equation*}
here it is automatically satisfied. We split the rest of the proof into five steps.

\vspace{2mm}

\noindent\textbf{Step I.} \emph{Convergence of $(Y^{T,n,x,a_0,a_\Gamma})_T$}. We begin by proving the following uniform estimate: $\P^{x,a_0,a_\Gamma}$-a.s.,
 	\begin{equation}\label{Sec:PDP_Sinfty_estimate}
 	0 \leq Y^{T,n,x,a_0,a_\Gamma}_s \leqslant   \frac{M_f}{\delta} +C^\ast \,M_c, \quad \forall \,s \in [0,\,T],
 	\end{equation}
 	 where $C^\ast$ is the constant defined in \eqref{Cdelta}. To this end, 
 for any $\bm\nu \in \mathcal{V}^n$ (the set of control maps $\bm\nu=(\nu^0,\nu^\Gamma)$, with both $\nu^0$ and $\nu^\Gamma$ bounded by $n$), let us introduce the compensated martingale measure $q^{\bm{\nu}}(ds \, dy \,db\,dc)= q(ds \, dy \,db\,dc)- [(\nu^0_s(b)-1)\,d_1(s,y,b,c)+(\nu^{\Gamma}_s(c)-1)\,d_2(s,y,b,c)]\,\tilde{p}(ds\,dy\,db\,dc)$ under $\P^{x,a_0,a_\Gamma}_{\bm\nu}$, with $d_1$ and $d_2$ given by respectively by \eqref{d1} and  \eqref{d2}.
 Taking  the expectation in \eqref{Sec:PDP_BSDE_penalized_TBIS_comp} under $\P^{x,a_0,a_\Gamma}_{\bm\nu}$,
 conditional to $\mathcal{F}_s$, and
 since   $Z^{T,n}$ is in  $\mathcal{G}^{\textbf{2},\beta}_{\textbf{x},\textbf{a}_0, \textbf{a}_\Gamma}(\textup{q};\textup{\textbf{0}},\textup{\textbf{T}})$, from Proposition \ref{Sec:PDP_lemma_P_nu_martingale} we get  that, $\P^{x,a_0,a_\Gamma}$-a.s.,
 \begin{align}
 &Y^{T,n}_{s}=- \spernuxa{\int_{(s,\,T]} \int_{A_0} e^{-\delta\,(r-s)}\,\{n [Z^{T,n}_{r}(X_r,b,J_r)]^- + \nu^0_r(b)\,Z^{T,n}_{r}(X_r,b,J_r)\}\, \lambda_0(db)\, dr\Big | \mathcal{F}_s}\nonumber \\
 & - \spernuxa{\int_{(s,\,T]} \int_{A_\Gamma} e^{-\delta\,(r-s)}\,\{n [Z^{T,n}_{r}(X_r,I_r,c)]^- + \nu^{\Gamma}_r(c)\,Z^{T,n}_{r}(X_r,I_r,c)\}\, \lambda_\Gamma(dc)\, dr\Big | \mathcal{F}_s}\notag\\ 
& + \spernuxa{\int_{(s,\,T]}
 	e^{-\delta\,(r-s)}\,f(X_{r},I_{r})\, dr+ \int_{(s,\,T]} e^{-\delta\,(r-s)}\,c(X_{r-},J_{r})\, d p^\ast_r\Big| \mathcal{F}_s}\,\, s \in [0,\,T]. \label{Sec:PDP_BSDE_nu}
 \end{align}
The right-hand side  of  estimate  \eqref{Sec:PDP_Sinfty_estimate} directly follows from 
 the elementary numerical inequality $n[z]^- + \nu z \geqslant 0$ for all $z \in \R$, $\nu \in (0,\,n]$, and the boundedness of $f$ and $c$.

Let us now prove that  $Y^{T,n}$ is nonnegative. To this end,  for $\varepsilon \in (0,\,1)$, let us consider the process ${\bm \nu}^\varepsilon := (\nu^{0,\varepsilon}, \nu^{\Gamma,\varepsilon}) \in {{\mathcal{V}}}^n$ defined by:
 \begin{align}
 \nu^{0,\varepsilon}_s(b) &= n\,\one_{\{Z_s^{T,n}(X_{s-},b,J_{s-}) \leqslant 0 \}} + \varepsilon\,\one_{\{0 < Z_s^{T,n}(X_{s-},b,J_{s-}) < 1\}} + \varepsilon\,Z_s^{T,n}(X_{s-},b,J_{s-})^{-1}\,\one_{\{Z_s^{n}(X_{s-},b,J_{s-}) \geqslant 1\}},\label{nuepsb}\\
  \nu^{\Gamma,\varepsilon}_s(c) &= n\,\one_{\{Z_s^{T,n}(X_{s-},I_{s-},c) \leqslant 0 \}} + \varepsilon\,\one_{\{0 < Z_s^{T,n}(X_{s-},I_{s-},c) < 1\}} + \varepsilon\,Z_s^{T,n}(X_{s-},I_{s-},c)^{-1}\,\one_{\{Z_s^{T,n}(X_{s-},I_{s-},c) \geqslant 1\}}.\label{nuepsc}
 \end{align}
 By construction, we have
 \begin{align*}
 n [Z_{s}^{T,n}(X_{s-},b, J_{s-})]^-+ \nu^{0,\varepsilon}_s(b)\,Z_{s}^{n}(X_{s-},b, J_{s-}) &\leqslant \varepsilon, \quad s \geqslant 0,\,b \in A_0,\\
  n [Z_{s}^{T,n}(X_{s-},I_{s-},c)]^-+ \nu^{\Gamma,\varepsilon}_s(c)\,Z_{s}^{n}(X_{s-},I_{s-},c) &\leqslant \varepsilon, \quad s \geqslant 0,\,c \in A_\Gamma.
 \end{align*}
Thus for the choice of $\bm \nu = {\bm \nu}^\varepsilon$ in \eqref{Sec:PDP_BSDE_nu}, we obtain
 \begin{align*}
 &Y_{s}^{T,n} \geqslant  -\varepsilon\, \frac{1-e^{-\delta (T-s)}}{\delta}\,(\lambda_\Gamma(A_\Gamma)+\lambda_0(A_0))\notag\\
  &+ \spernuxaepsilon{ \int_{(s,\,T]}
 	e^{-\delta\,(r-s)}\,f(X_{r},I_{r})\, dr+ \int_{(s,\,T]} e^{-\delta\,(r-s)}\,c(X_{r-},J_{r})\, d p^\ast_r \Big | \mathcal{F}_s}.
 \end{align*}
Since $f,c$ are positive, it follows that
 \begin{align}
 	Y_{s}^{T,n} 
 	&\geqslant \essinf_{\bm \nu \in \mathcal{V}^n} \spernuxa{\int_{(s,\infty)} e^{-\delta\,(r-s)}\,f(X_{r},I_{r})\, dr +\int_{(s,\infty)} e^{-\delta\,(r-s)}\,c(X_{r-},J_{r-})\, d p^\ast_r\Big | \mathcal{F}_s}\notag\\
 	&-\frac{\varepsilon}{\delta}\,(\lambda_0(A_0)+\lambda_\Gamma(A_\Gamma)).\label{secondineq}
 \end{align}
We conclude by the arbitrariness of $\varepsilon$.

Now, let us study the convergence of $(Y^{T,n})_T$. 
Take $T,T'>0$, with $T< T'$, and $s \in [0,\,T]$. Then
\begin{eqnarray}\label{Sec:PDP_Ybar_conv}
 |Y_s^{T',n}-Y_s^{T,n}|^2
 \leqslant e^{-2\,\delta\,(T-s)}\,\spernuxaepsilon{|Y_T^{T',n}-Y_T^{T,n}|^2|\mathcal{F}_s}\overset{T,\,T' \rightarrow \infty}{\longrightarrow} 0,
 \end{eqnarray}
 where the convergence result follows from \eqref{Sec:PDP_Sinfty_estimate}.
 Let us now consider the sequence of real-valued c\`adl\`ag adapted processes $(Y^{T,n})_T$. It follows from \eqref{Sec:PDP_Ybar_conv} that, for any $t\geqslant 0$, the sequence $(Y_t^{T,n}(\omega))_T$ is  Cauchy  for almost every $\omega$, so that it converges $\P^{x,a}$-a.s. to some $\mathcal{F}_t$-measurable random variable $Y_t^n$, which is bounded by the right-hand side of \eqref{Sec:PDP_Sinfty_estimate}. 
 Moreover, using again \eqref{Sec:PDP_Ybar_conv} and \eqref{Sec:PDP_Sinfty_estimate}, we see that, for any $0\leqslant S< T\wedge T'$, with $T,T' >0$, we have
 \begin{equation}\label{Sec:PDP_conv_Y}
 \sup_{0\leqslant t\leqslant S}|Y_t^{T',n}-Y_t^{T,n}| \leqslant e^{-\delta\,(T\wedge T' -S)}\,\left(\frac{M_f}{\delta} + C^\ast\,M_c\right)\overset{T,T'\rightarrow \infty}{\longrightarrow} 0.
 \end{equation}
 Since each $Y^{T,n}$ is a c\`adl\`ag process, it follows that $Y^n$ is c\`adl\`ag, as well. Finally, from estimate \eqref{Sec:PDP_Sinfty_estimate} we see that  $Y^n$ is uniformly bounded and therefore belongs to $\textup{\textbf{S}}^\infty$.

\vspace{2mm}
  
\noindent\textbf{Step II.} \emph{Convergence of $(Z^{T,n,x,a_0,a_\Gamma})_T$.} 
 Let $S,\,T,\,T'>0$, with $S<T<T'$.  Then, applying It\'o's formula to $e^{-2\,\delta\,t}|Y^{T',n}_t-Y^{T,n}_t|^2$ between $0$ and $S$,  and taking the expectation, we get
\begin{align*}
&\frac{1}{2}\sperxa{\sum_{r \in (0,\,S]}e^{-2\,\delta\,r}
\left|\int_{E \times A_0\times A_\Gamma} (Z^{T',n}_r(y,b,c)-Z^{T,n}_r(y,b,c))\,q(\{r\}\times dy \,db\,dc)\right|^2}\\
 &  \leqslant e^{-2\,\delta\,S}\sperxa{|Y^{T',n}_S-Y^{T,n}_S|^2} \\
 & +4(n^2+1)\,(\lambda_0(A_0)+\lambda_\Gamma(A_\Gamma)) \,\sperxa{\int_0^S e^{-2\,\delta\,r}\,|Y^{T',n}_r-Y^{T,n}_r|^2\,dr}\overset{T,T'\rightarrow \infty}{\longrightarrow}0,
 \end{align*}
 where the convergence to zero follows from estimate \eqref{Sec:PDP_conv_Y}.
 Then, for any $S>0$, we see  that $(Z^{T,n}_{|[0,\,S]})_{T>S}$ is a Cauchy sequence in the Hilbert space $\mathcal{G}^{\textbf{2}}_{\textbf{x},\textbf{a}_0, \textbf{a}_\Gamma}(\textup{q};\textbf{0},\,\textbf{S})$. Therefore, we deduce the existence of $Z^n\in\mathcal{G}^{\textbf{2}}_{\textbf{x},\textbf{a}_0,\textbf{a}_\Gamma,\textbf{\textup{loc}}}(\textup{q})$ such that $(Z^{T,n}_{|[0,\,S]})_{T>S}$ converges to $Z^n_{|[0,\,S]}$ in $\mathcal{G}^{\textbf{2}}_{\textbf{x},\textbf{a}_0, \textbf{a}_\Gamma}(\textup{q};\textbf{0},\,\textbf{S})$. Hence, from the convergence of $(Y^{T,n})_T$ and $(Z^{T,n})_n$, we can pass to the limit in equation \eqref{Sec:PDP_BSDE_penalized_T} as $T \rightarrow \infty$, from which we deduce that $(Y^n,Z^n)$ (also denoted as $(Y^{n,x,a_0,a_\Gamma},Z^{n,x,a_0,a_\Gamma})$) solves the following penalized BSDE on infinite horizon: $\P^{x,a_0,a_\Gamma}$-a.s.,
 \begin{align}\label{Sec:PDP_BSDE_penalized}
 Y_{s}^{n,x,a_0,a_\Gamma}&= Y_{T}^{n,x,a_0,a_\Gamma}-\delta\int_{(s,T]} Y_r^{n,x,a_0,a_\Gamma}\,dr + \int_{(s,\,T]} f(X_{r},I_{r})\, dr+ \int_{(s,\,T]} c(X_{r-},J_{r-})\, d p^\ast_r\nonumber\\
 & - \int_{(s,\,T]} \int_{A_0} Z_{r}^{n,x,a_0,a_\Gamma}(X_{r},b,J_r) \,\lambda_{0}(db)\, dr - \int_{(s,\,T]} \int_{A_\Gamma} Z_{r}^{n,x,a_0,a_\Gamma}(X_{r},I_r,c) \,\lambda_{\Gamma}(dc)\, dr \nonumber\\
 & - \big(K_T^{n,x,a_0,a_\Gamma} - K_s^{n,x,a_0,a_\Gamma}\big) - \int_{(s,\,T]} \int_{E \times A_0\times A_\Gamma} Z_{r}^{n,x,a_0,a_\Gamma}(y,b,c) \, q(dr\,dy\,db\,dc),
\end{align}
for all $0\leqslant s\leqslant T<\infty$, where
\[
 K_s^{n,x,a_0,a_\Gamma} := n\int_0^s \left(\int_{A_0}[Z_{r}^{n,x,a_0,a_\Gamma}(X_r,b,J_r)]^{-}\,\lambda_0(db)+\int_{A_\Gamma}[Z_{r}^{n,x,a_0,a_\Gamma}(X_r,I_r,c)]^{-}\,\lambda_\Gamma(dc)\right)dr.
\]
Notice that equation \eqref{Sec:PDP_BSDE_penalized} can also be written as follows:
\begin{align}\label{Sec:PDP_BSDE_penalized_TBIS_comp_n}
 	 Y_{s}^{n,x,a_0,a_\Gamma}
 	&=Y_{T}^{n,x,a_0,a_\Gamma}\,e^{-\delta (T-s)} + \int_{(s,\,T]} \tilde{f}^n(r-s,X_{r-},I_{r-},\,J_{r-},\, Z_r^{n,x,a_0,a_\Gamma})\, d A_r\nonumber\\ 
 	&- \int_{(s,\,T]} e^{-\delta (r-s)}\int_{E \times A_0\times A_{\Gamma}}  Z_{r}^{n,x,a_0,a_\Gamma}(y,b,c) \, q(dr\,dy\,db\,dc),\quad s \in [0,\,T],
 	\end{align}
where	   
	$\tilde{f}^n$ is the deterministic function defined in \eqref{Sec:PDP_tildef^n}.
 
\vspace{2mm}

\noindent\textbf{Step III.} \emph{Representation formula for $Y^{n,x,a_0,a_\Gamma}$.} Our aim is to prove the following representation formula:
\begin{equation}\label{Sec:PDP_rep_Y_n}
Y_{s}^{n,x,a_0,a_\Gamma} = \essinf_{\bm \nu \in \mathcal{V}^n} \spernuxa{\int_{(s,\infty)} e^{-\delta\,(r-s)}\,f(X_{r},I_{r})\, dr +\int_{(s,\infty)} e^{-\delta\,(r-s)}\,c(X_{r-},J_{r-})\, d p^\ast_r\Big | \mathcal{F}_s},
\end{equation}
for all $s\geq0$. As at the beginning of \textbf{Step I}, for any $\bm\nu \in \mathcal{V}^n$, we consider the compensated martingale measure $q^{\bm{\nu}}(ds \, dy \,db\,dc)= q(ds \, dy \,db\,dc)- [(\nu^0_s(b)-1)\,d_1(s,y,b,c)+(\nu^{\Gamma}_s(c)-1)\,d_2(s,y,b,c)]\,\tilde{p}(ds\,dy\,db\,dc)$ under $\P^{x,a_0,a_\Gamma}_{\bm\nu}$. We take  the expectation in \eqref{Sec:PDP_BSDE_penalized_TBIS_comp_n} under $\P^{x,a_0,a_\Gamma}_{\bm\nu}$,
 conditional to $\mathcal{F}_s$. 
For every $T>0$, recalling that   $Z^{n}$ is in  $\mathcal{G}^{\textbf{2}}_{\textbf{x},\textbf{a}_0, \textbf{a}_\Gamma}(\textup{q};\textup{\textbf{0}},\textup{\textbf{T}})$, 
from Proposition \ref{Sec:PDP_lemma_P_nu_martingale}
we get
 \begin{align}
 Y_{s}^{n}= & - \spernuxa{\int_{(s,\,T]} \int_{A_0} e^{-\delta\,(r-s)}\,\{n [Z^{n}_{r}(X_r,b,J_r)]^- + \nu^0_r(b)\,Z^{n}_{r}(X_r,b,J_r)\}\, \lambda_0(db)\, dr\Big | \mathcal{F}_s}\label{Sec:PDP_BSDE_spernu} \\
 & - \spernuxa{\int_{(s,\,T]} \int_{A_\Gamma} e^{-\delta\,(r-s)}\,\{n [Z^{n}_{r}(X_r,I_r,c)]^- + \nu^{\Gamma}_r(c)\,Z^{n}_{r}(X_r,I_r,c)\}\, \lambda_\Gamma(dc)\, dr\Big | \mathcal{F}_s}\notag\\
  & +\spernuxa{e^{-\delta\,(T-s)}\,Y_{T}^{n,x,a_0,a_\Gamma} + \int_{(s,\,T]}
 	e^{-\delta\,(r-s)}\,f(X_{r},I_{r})\, dr+ \int_{(s,\,T]} e^{-\delta\,(r-s)}\,c(X_{r-},J_{r})\, d p^\ast_r \Big | \mathcal{F}_s}.\notag
 \end{align}
 From the elementary inequality $n[z]^- + \nu z \geqslant 0$, $z \in \R$, $\nu \in (0,\,n]$, we deduce
 \begin{align*}
 	Y_s^{n} 
 	\leq 
 	 \spernuxa{e^{-\delta\,(T-s)}\,Y_{T}^{n} + \int_{(s,\,T]}
 	e^{-\delta\,(r-s)}\,f(X_{r},I_{r})\, dr+ \int_{(s,\,T]} e^{-\delta\,(r-s)}\,c(X_{r-},J_{r})\, d p^\ast_r \Big | \mathcal{F}_s}.
 \end{align*}
 Since $Y^{n}$ is in $\textbf{S}^{\infty}$,
  sending $T \rightarrow \infty$, we obtain, by the conditional version of the Lebesgue dominated convergence theorem,
 \begin{align*}
 		Y_s^{n} &\leq \spernuxa{
 	\int_{(s,\,\infty)}
 	e^{-\delta\,(r-s)}\,f(X_{r},I_{r})\, dr+ \int_{(s,\,\infty)} e^{-\delta\,(r-s)}\,c(X_{r-},J_{r})\, d p^\ast_r \Big | \mathcal{F}_s}.
 \end{align*}
Hence 
 \begin{align}
 Y_s^{n} &\leq  \essinf_{{\bm \nu} \in \mathcal{V}^n}\spernuxa{
 	\int_{(s,\,\infty)}
 	e^{-\delta\,(r-s)}\,f(X_{r},I_{r})\, dr+ \int_{(s,\,\infty)} e^{-\delta\,(r-s)}\,c(X_{r-},J_{r})\, d p^\ast_r \Big | \mathcal{F}_s}.
 \label{Sec:PDP_ineq_esssup}
 \end{align}
 On the other hand, for $\varepsilon \in (0,\,1)$, let us consider the process ${\bm \nu}^\varepsilon := (\nu^{0,\varepsilon}, \nu^{\Gamma,\varepsilon}) \in {{\mathcal{V}}}^n$ defined by  \eqref{nuepsc}-\eqref{nuepsc}, with $Z^{T,n}$ replaced by $Z^{n}$.
Thus for this choice of $\bm \nu = {\bm \nu}^\varepsilon$ in \eqref{Sec:PDP_BSDE_spernu}, we obtain
 \begin{align*}
 &Y_{s}^{n} \geqslant  -\varepsilon\, \frac{1-e^{-\delta (T-s)}}{\delta}\,(\lambda_\Gamma(A_\Gamma)+\lambda_0(A_0))\notag\\
  &+ \spernuxaepsilon{e^{-\delta\,(T-s)}\,Y_{T}^{n} + \int_{(s,\,T]}
 	e^{-\delta\,(r-s)}\,f(X_{r},I_{r})\, dr+ \int_{(s,\,T]} e^{-\delta\,(r-s)}\,c(X_{r-},J_{r})\, d p^\ast_r \Big | \mathcal{F}_s}.
 \end{align*}
 Letting $T \rightarrow \infty$, since $f,c$ are bounded  and  $Y^{n,x,a_0,a_\Gamma} \in\textbf{S}^{\infty}$, it follows that
 \begin{align}
 	Y_{s}^{n} 
 	&\geqslant \essinf_{\bm \nu \in \mathcal{V}^n} \spernuxa{\int_{(s,\infty)} e^{-\delta\,(r-s)}\,f(X_{r},I_{r})\, dr +\int_{(s,\infty)} e^{-\delta\,(r-s)}\,c(X_{r-},J_{r-})\, d p^\ast_r\Big | \mathcal{F}_s}\notag\\
 	&-\frac{\varepsilon}{\delta}\,(\lambda_0(A_0)+\lambda_\Gamma(A_\Gamma)).\label{secondineq}
 \end{align}
 Taking into account the arbitrariness of $\varepsilon$, 
 the required representation of $Y^{n,x,a_0,a_\Gamma}$ follows from \eqref{Sec:PDP_ineq_esssup} and \eqref{secondineq}.
 
\vspace{2mm}

\noindent\textbf{Step IV.} \emph{Uniform estimate on $(Z^{n,x,a_0,a_\Gamma},K^{n,x,a_0,a_\Gamma})_n$.} Let us prove that, for every $T >0$, there exists a constant $C$, depending only on $M_f$, $M_c$, $\delta$,   $T$, $C^\ast$, such that
\begin{eqnarray}\label{Sec:PDP_estimate_bsde}
||Z^{n,x,a_0,a_\Gamma}||^2_{\mathcal{G}^{\textbf{2}}_{\textbf{x},\textbf{a}_0, \textbf{a}_\Gamma}(\textup{q};\textup{\textbf{0}},\,\textup{\textbf{T}})} + ||K^{n,x,a_0,a_\Gamma}||^2_{\textup{\textbf{K}}^{\textbf{2}}_{\textbf{x},\textbf{a}_0, \textbf{a}_\Gamma}(\textup{\textbf{0}},\,\textup{\textbf{T}})} \leqslant C.
\end{eqnarray}
Fix $T >0$. In what follows we shall denote by $C >0$ a generic positive constant depending on $M_f$, $M_c$, $C^\ast$,  $\delta$ and  $T$,
 which may vary from line to line.
To simplify notation, we denote $Y^{n,x,a_0,a_\Gamma}$, $Z^{n,x,a_0,a_\Gamma}$, $K^{n,x,a_0,a_\Gamma}$ simply by $Y^n$, $Z^n$, $K^n$. Applying It\^o's formula to $|Y_s^n|^2$ between $0$ and $T$, and taking the expectation with respect to $\P^{x,a_0,a_{\Gamma}}$, recalling also Lemma \ref{L_simplifiedNorm}, we obtain 
 \begin{align*}
&\sperxa{\int_{(0,T]}
 \int_{E \times A_0 \times A_\Gamma}\ \big|Z^{n}_s(y,b,c) - \hat Z^{n}_s \,\one_{K}(s)\big|^2\,\tilde p(ds\,dy\,db)}\\
 &\leq \sperxa{
 |Y_T^{n}|^2}  -2 \sperxa{\int_{(0,T]} 
 Y_{s}^{n}\,dK^{n}_s} 
\notag\\	
&+2 \sperxa{\sum_{s \in (0,\,T]}
\left(\int_{E \times A_0\times A_\Gamma} Z_s^{n}(y,b,c)\,q(\{s\}\times dy \,db\,dc)\right)c(X_{s-},J_{s-})\,\one_{X_{s-} \in \partial E}	}\notag\\
 & +2\sperxa{\int_{(0,\,T]}
 Y_{s}^{n}\,f(X_s,I_s)\,ds}+ 2\,\sperxa{\int_{(0,\,T]}
 Y^{n}_{s-}\,c(X_{s-},J_{s-})\, d p^\ast_s}.
 \end{align*}
Using the elementary inequality $2\,a\,b \leq \gamma\, a^2 + \frac{1}{\gamma}\,b^2$, with 
 $\gamma \in \R_+ \setminus \{0\}$, $\gamma <1$, we get
  \begin{align*}
&(1-\gamma)\,\sperxa{\int_{(0,T]}
 \int_{E \times A_0 \times A_\Gamma}\ \big|Z^{n}_s(y,b,c) - \hat Z^{n}_s \,\one_{K}(s)\big|^2\,\tilde p(ds\,dy\,db)}\\
 &\leq \sperxa{
 |Y_T^{n}|^2}  -2 \sperxa{\int_{(0,T]} 
 Y_{s}^{n}\,dK^{n}_s} 
+\frac{1}{\gamma} \,\sperxa{\sum_{s \in (0,\,T]}
|c(X_{s-},J_{s-})|^2\,\one_{X_{s-} \in \partial E}	}\notag\\
 & +2\sperxa{\int_{(0,\,T]}
 Y_{s}^{n}\,f(X_s,I_s)\,ds}+ 2\,\sperxa{\int_{(0,\,T]}
 Y^{n}_{s-}\,c(X_{s-},J_{s-})\, d p^\ast_s}.
 \end{align*}
  Set  now 
 $C_Y := \frac{M_f}{\delta} + C^\ast \,M_c$. Recalling the uniform estimate \eqref{Sec:PDP_Sinfty_estimate}   on $Y^n$, we obtain
 \begin{align}\label{Sec:PDP_partial_estimate_Z}
 &(1-\gamma)\,\sperxa{\int_{(0,T]}
 \int_{E \times A_0 \times A_\Gamma}\ \big|Z^{n}_s(y,b,c) - \hat Z^{n}_s \,\one_{K}(s)\big|^2\,\tilde p(ds\,dy\,db\,dc)}  \nonumber\\
 &\leqslant 
 \frac{1}{\gamma}\,\,M_c^2 \,C^\ast(T)
  + 
 C_Y^2 
 \, +2\,C_Y\,(M_f \,T +M_c\,C^\ast(T))
 + 2\,C_Y\,\sperxa{ K_T^{n}},
 \end{align}
 where $C^\ast(t)$ is the deterministic function defined in \eqref{conv_exp_sper_pstar}.
On the other hand, from  \eqref{Sec:PDP_BSDE_penalized}, we get
 \begin{align}\label{eq_Kn}
  K_T^{n}  &= Y_T^{n} - Y_0^{n} -\delta\, \int_{(0,\,T]}Y_s^{n,x,a}\,ds + \int_{(0,\,T]} f(X_s,I_s)\,ds + \int_{(0,\,T]} c(X_{s-},J_{s-})\,d p^\ast_s
\nonumber\\
 &  - \int_{(0,\,T]}\int_{A_0} \,Z_s^{n}(X_s,b, J_s)\,\lambda_0(db)\,ds -  \int_{(0,\,T]}\int_{A_\Gamma} \,Z_s^{n}(X_s, I_s,c)\,\lambda_0(dc)\,ds \notag\\
 &- \int_{(0,\,T]}\int_{E \times A_0 \times A_\Gamma} \,Z_s^{n}(y,b,c)\,q(ds\,dy\,db\,dc).
 \end{align}
Using again the inequality  $2ab \leqslant \frac{1}{\eta}a^2 + \eta b^2$, for any $\eta=\alpha, k 
 >0$, and taking the expectation in \eqref{eq_Kn}, we find 
 \begin{align}\label{Sec:PDP_partial_estim_K^n}
 &2\,\sperxa{K_T^{n}}  \leqslant 4\,C_Y +  2\,\delta\,C_Y\,T
 + 2\,M_f\,T + 2\,M_c\,C^\ast(T)\notag\\
 &+ \frac{T}{\alpha}\,\lambda_0(A_0)+ \alpha\,\sperxa{\int_{(0,\,T]}\int_{A_0} \,|Z_s^{n}(X_s,b,J_s)|^2\,\lambda_0(db)\,ds}\notag\\
  &+ \frac{T}{k}\,\lambda_\Gamma(A_\Gamma)+ k\,\sperxa{\int_{(0,\,T]}\int_{A_\Gamma} \,|Z_s^{n}(X_s,I_s, c)|^2\,\lambda_\Gamma(dc)\,ds}.
 \end{align}
 Plugging \eqref{Sec:PDP_partial_estim_K^n} into \eqref{Sec:PDP_partial_estimate_Z},
 we obtain
 \begin{align*}
& (1-\gamma)\,\sperxa{\int_{(0,T]}
 \int_{E \times A_0 \times A_\Gamma}\ \big|Z^{n}_s(y,b,c) - \hat Z^{n}_s \,\one_{K}(s)\big|^2\,\tilde p(ds\,dy\,db)}
 \leqslant C+ \\
 &+(\alpha \vee k)\,C_Y\left(1+ 2\,T\right)\left(\int_{(0,\,T]}\left[\int_{A_0} \,|Z_s^{n}(X_s,b,J_s)|^2\,\lambda_0(db) + \int_{A_\Gamma} \,|Z_s^{n,x,a}(X_s,I_s,c)|^2\,\lambda_\Gamma(dc)\right]ds\right).
 \end{align*}
Choosing $\alpha = k =\frac{1-\gamma}{2\,C_Y\,(1+2\,T)}$, we get
 \begin{eqnarray*}
 	 (1-\gamma)\,\sperxa{\int_{(0,T]}
 \int_{E \times A_0 \times A_\Gamma}\ \big|Z^{n}_s(y,b,c) - \hat Z^{n}_s \,\one_{K}(s)\big|^2\,\tilde p(ds\,dy\,db)}
 	\leqslant
 	C,
 \end{eqnarray*}
which gives the required uniform estimate for $(Z^{n})_n$,
 and also for $(K^{n})_n$ by \eqref{eq_Kn}.

\vspace{2mm}

\noindent\textbf{Step V.} \emph{Convergence of $(Y^{n,x,a_0,a_\Gamma},Z^{n,x,a_0,a_\Gamma},K^{n,x,a_0,a_\Gamma})_n$.} It follows from estimate \eqref{Sec:PDP_Sinfty_estimate} and the representation formula \eqref{Sec:PDP_rep_Y_n}, that the sequence $(Y^n)_n$ converges in a nondecreasing way to some uniformly bounded process $Y$. By \eqref{Sec:PDP_rep_Y_n}, we then deduce the representation formula \eqref{Sec:PDP_rep_Y} for $Y$. In addition, by the uniform estimate \eqref{Sec:PDP_estimate_bsde} it follows that there exist $Z^{x,a_0,a_\Gamma}\in\mathcal{G}^{\textbf{2}}_{\textbf{x},\textbf{a}_0,\textbf{a}_\Gamma,\textbf{\textup{loc}}}(\textup{q}) $ and a nondecreasing, predictable process $K^{x,a_0,a_\Gamma}$, with $K_0=0$ and $\E^{x,a_0,a_\Gamma}[|K_T^{x,a_0,a_\Gamma}|^2]<\infty$, such that:
\begin{itemize}
\item $Z^{x,a_0,a_\Gamma}$ is the weak limit of $(Z^{n,x,a_0,a_\Gamma})_n$ in $\mathcal{G}^{\textbf{2}}_{\textbf{x},\textbf{a}_0,\textbf{a}_\Gamma,\textbf{\textup{loc}}}(\textup{q})$;
\item $K_{s}^{x,a_0,a_\Gamma}$ is the weak limit of $(K_{s}^{n,x,a_0,a_\Gamma})_n$ in $\textup{\textbf{L}}^{\textbf{2}}_{\textbf{x},\textbf{a}_0,\textbf{a}_\Gamma}(\mathcal{F}_{s})$, for every $s \geqslant 0$.
\end{itemize}
By Lemma 2.2 in \cite{Pe}, we deduce that both $Y^{x,a_0,a_\Gamma}$ and $K^{x,a_0,a_\Gamma}$ are c\`adl\`ag processes, so that $Y^{x,a_0,a_\Gamma}\in\textup{\textbf{S}}^\infty$ and $K^{x,a_0,a_\Gamma}\in\textup{\textbf{K}}^{\textbf{2}}_{\textbf{x},\textbf{a}_0,\textbf{a}_\Gamma,\textbf{\textup{loc}}}$. Letting $n\rightarrow\infty$ in equation \eqref{Sec:PDP_BSDE_penalized}, we see that $(Y^{x,a_0,a_\Gamma},Z^{x,a_0,a_\Gamma},K^{x,a_0,a_\Gamma})$ solves equation \eqref{Sec:PDP_BSDE}.

Consider now another solution $(\tilde{Y}, \tilde{Z},\tilde{K})\in \textbf{S}^{\infty}\times \mathcal{G}^{\textbf{2}}_{\textbf{x},\textbf{a}_0,\textbf{a}_\Gamma,\textbf{\text{loc}}}(\textup{q})\times \textbf{K}^{\textbf{2}}_{\textbf{x},\textbf{a}_0,\textbf{a}_\Gamma,\textbf{\text{loc}}}$ to \eqref{Sec:PDP_BSDE}-\eqref{Sec:PDP_BSDE_constraint1}-\eqref{Sec:PDP_BSDE_constraint2}. Then, it is quite easy to check that
\[
\tilde Y_s^{x,a_0,a_\Gamma} \leq \essinf_{\bm\nu \in \mathcal{V}}\spernuxa{\int_{(s,\infty)} e^{-\delta (r-s)}\,f(X_{r},I_{r})\, dr +\int_{(s,\infty)} e^{-\delta (r-s)}\,c(X_{r-},J_{r-})\, d p^\ast_r\Big| \mathcal{F}_s},
\]
for all $s\geq0$. This implies the maximality of $(Y^{x,a_0,a_\Gamma},Z^{x,a_0,a_\Gamma},K^{x,a_0,a_\Gamma})$.

Concerning the jump constraints, we simply notice that they  are a direct consequence of the uniform estimate \eqref{Sec:PDP_estimate_bsde} on the norm $||K^{n,x,a_0,a_\Gamma}||^2_{\textup{\textbf{K}}^{\textbf{2}}_{\textbf{x},\textbf{a}_0, \textbf{a}_\Gamma}(\textup{\textbf{0}},\,\textup{\textbf{T}})}$.

Finally, regarding the uniqueness result, 
 let $(Y,Z,K)$ and $(Y',Z',K')$ be two maximal solutions of \eqref{Sec:PDP_BSDE}-\eqref{Sec:PDP_BSDE_constraint1}-\eqref{Sec:PDP_BSDE_constraint2}. 
 The component $Y$ is unique by definition.
 Let us now consider the difference between the two backward equations. We get: $\P^{x,a_0,a_\Gamma}$-a.s.
 \begin{align}\label{Sec:PDP_rewritten}
 	&\int_{(0,\,t]} \int_{E \times A_0\times A_\Gamma} (Z_{s}(y,b,c)- Z'_{s}(y,b,c))\, q(ds\,dy\,db\,dc) \notag \\
 	&= (K_t -K'_t) -\int_{(0,\,t]} \int_{A_0} (Z_{s}(X_s,b,J_s)- Z'_{s}(X_s,b,J_s))\, \lambda_0(db)\,ds\\
 	&-\int_{(0,\,t]} \int_{A_\Gamma} (Z_{s}(X_s,I_s,c)- Z'_{s}(X_s,I_s,c))\, \lambda_\Gamma(dc)\,ds, \quad 0\leqslant t\leqslant T<\infty.\notag
 \end{align}
 The right-hand side of \eqref{Sec:PDP_rewritten} is a predictable process, therefore it has no totally inaccessible jumps (see, e.g., Proposition 2.24, Chapter I, in \cite{JS}); on the other hand, by Lemma \ref{P_simplifiedNorm},  together with \eqref{D}-\eqref{K}, it follows that the left-hand side of \eqref{Sec:PDP_rewritten} is a jump process with only totally inaccessible jumps.
 This implies that $Z=Z'$ in $\mathcal{G}^{\textbf{2}}_{\textbf{x},\textbf{a}_0,\textbf{a}_\Gamma,\textbf{\text{loc}}}(\textup{q})$, and as a consequence the component $K$ is unique as well.
\endproof

\section{A BSDE representation for the value function}
\label{Sec:PDP_Section_nonlinear_IPDE}

The aim of the present section is to prove that the value function $V$ in \eqref{Sec:PDP_value_function} can be represented in terms of the maximal solution to the BSDE with nonnegative jumps \eqref{Sec:PDP_BSDE}-\eqref{Sec:PDP_BSDE_constraint1}-\eqref{Sec:PDP_BSDE_constraint2}. Firstly, we introduce the deterministic function $v: E \times A_0 \times A_\Gamma \rightarrow \R$ given by
\begin{equation}
v(x,a_0,a_\Gamma):= Y_0^{x,a_0,a_\Gamma}, \quad (x,a_0,a_\Gamma) \in E\times A_0 \times A_\Gamma. \label{Sec:PDP_def_v}
\end{equation}
\begin{theorem}\label{Sec:PDP_THm_Feynman_Kac_HJB}
Let assumptions \textup{\textbf{(H$\textup{h$\lambda$QR}$)}},  {\bf (H$0$)}, {\bf (H$0'$)}, \textup{\textbf{(H$\lambda_0\lambda_\Gamma$)}} and \textup{\textbf{(H$\textup{fc}$)}} hold. Then, the function $v$ in \eqref{Sec:PDP_def_v} does not depend on its last arguments:  
\begin{equation}
 v(x,a_0,a_\Gamma)= v(x,a'_0,a'_\Gamma), \quad x \in E,\,\, (a_0,a'_0)\in A_0,\,\,(a_\Gamma,a'_\Gamma)\in A_\Gamma.\label{Sec:PDP_vdelta_not_dep_a}
\end{equation}
By an abuse of notation, we define  the function $v$ on $E$ by
 \begin{equation}\label{Sec:PDP_redefine_vdelta}
 v(\cdot):= v(\cdot,a_0,a_\Gamma), \quad 	\textup{for any}\,\, (a_0,a_\Gamma) \in A_0 \times A_\Gamma.
 \end{equation}
Then $v$ is continuous and bounded. Moreover, $v$ admits the representation formula: $\P^{x,a_0,a_\Gamma}$-a.s.,
\begin{equation}\label{Sec:PDP_ident_vdelta}
v(X_s) = Y_s^{x,a_0,a_\Gamma}, \quad \forall s \geq 0.
\end{equation}
\end{theorem}
\proof
We split the proof into three steps.

\vspace{2mm}

\noindent\textbf{Step I.} \emph{The identification property of  $Y^{x,a_0,a_\Gamma}$.} A first fundamental preliminary result we have to prove is the following identification property: 
 	for every $(x,a_0,a_\Gamma)\in E \times A_0 \times A_\Gamma$, 
 	$\P^{x,a_0,a_\Gamma}$-a.s.,	$s \geq 0$, 
 	\begin{equation}
 	Y_s^{x,a_0,a_\Gamma}=v(X_s, I_s, J_s), \label{Sec:PDP_ident_vdelta1}
 	\end{equation}
where $v$ is the deterministic function defined by \eqref{Sec:PDP_def_v}. Recall from the proof of Theorem \ref{Sec:PDP_Thm_ex_uniq_maximal_BSDE} that $Y^{x,a_0,a_\Gamma}$ is constructed from $Y^{T,n,x,a_0,a_\Gamma}$ (see equation \eqref{Sec:PDP_BSDE_penalized_T}), taking firstly the limit as $T\rightarrow\infty$, and then as $n\rightarrow\infty$. Therefore, it is enough to prove property \eqref{Sec:PDP_ident_vdelta1} for $Y^{T,n,x,a_0,a_\Gamma}$. For simplicity of notation, denote the pair $(Y^{T,n,x,a_0,a_\Gamma},Z^{T,n,x,a_0,a_\Gamma})$, solution to equation \eqref{Sec:PDP_BSDE_penalized_T}, simply as $(Y^{T,n},Z^{T,n})$. Then, we know from the fixed point argument 
 giving the well-posedness of the penalized BSDE \eqref{Sec:PDP_BSDE_penalized} (see the proof of Theorem 4.1 in \cite{BandiniBSDE})
   that there exists a sequence $(Y^{T,n,k},Z^{T,n,k})_k$
   in $\textbf{S}^{\infty}
   \times \mathcal{G}^{\textbf{2}}_{{\textbf{x},\textbf{a}_0, \textbf{a}_\Gamma}}(\textup{q}; \textbf{0,\,S})$ converging to $(Y^{T,n},Z^{T,n})$ in $\textbf{S}^{\infty}
   \times \mathcal{G}^{\textbf{2}}_{{\textbf{x},\textbf{a}_0, \textbf{a}_\Gamma}}(\textup{q}; \textbf{0,\,S})$, such that $(Y^{T,n,0}, Z^{T,n,0})= (0,0)$ and
 \begin{align}\label{kBSDE}
 	Y_{t}^{T,n,k+1}&= Y_{S}^{T,n,k}-\delta\int_t^S Y_r^{T,n,k}\,dr + \int_{(t,\,S]} f(X_{r},I_{r})\, dr + \int_{(t,\,S]} c(X_{r-},J_{r-})\, d p^\ast_r\nonumber\\
 	& -n\int_t^S\int_{A_0} [Z_r^{T,n,k}(X_r,b,I_r)]^-\,\lambda_0(db)\,dr - \int_{(t,\,S]} \int_{A_0} Z_{r}^{T,n,k}(X_{r},b,I_r) \,\lambda_{0}(db)\, dr, \nonumber\\
 	& -n\int_{(t,\,S]}\int_{A_\Gamma} [Z_r^{T,n,k}(X_r,I_r,c)]^-\,\lambda_\Gamma(dc)\,dr - \int_{(t,\,S]} \int_{A_\Gamma} Z_{r}^{T,n,k}(X_{r},I_r,c) \,\lambda_{\Gamma}(dc)\, dr\nonumber\\
 	& - \int_{(t,\,S]} \int_{E \times A_0 \times A_\Gamma} Z_{r}^{T,n,k+1}(y,b,c) \, q(dr\,dy\,db\,dc),\quad \quad 0 \leqslant t\leqslant S.
 \end{align}
 Let us define
 $v^{T,n}(x,a_0,a_\Gamma):= Y_{0}^{T,n}$, $v^{T,n,k}(x,a_0,a_\Gamma):= Y_{0}^{T,n,k}$.
 For $k=0$, we have, $\P^{x,a_0,a_\Gamma}$-a.s.,
 $$
 Y_{t}^{T,n,1}= \sperxa{\int_{(t,\,S]} f(X_{r},I_{r})\, dr + \int_{(t,\,S]} c(X_{r-},J_{r-})\, d p^\ast_r \Big |\mathcal{F}_t}, \quad t \in [0,\,S].
 $$
 Then, from the Markov property of $(X,I,J)$ we get, $\P^{x,a_0,a_\Gamma}$-a.s.,
 $
 Y_{t}^{T,n,1}= v^{T,n,1}(X_{t},I_{t},J_{t})$,
 and in particular
  \begin{align*}
 	\Delta Y_{t}^{T,n,1}&=  -c(X_{t-},J_{t-})\, \Delta p^\ast_t
 	+ \int_{E \times A_0 \times A_\Gamma} Z_{t}^{T,n,1}(y,b,c) \, q(\{t\}\times \,dy\,db\,dc)\\
 	&=  -c(X_{t-},J_{t-})\, \Delta p^\ast_t
 	+ Z_{t}^{T,n,1}(X_t,I_t, J_t)- \hat Z^{T,n,1}_t\, \Delta p^\ast_t,\quad \quad 0 \leqslant t\leqslant S.
 \end{align*}
 This gives
 \begin{equation*}
 Z_{t}^{T,n,1}(y,b,c)-\hat Z_{t}^{T,n,1}\,\one_{X_{t-} \in \partial E} = v^{T,n,1}(y,b,c)-v^{T,n,1}(X_{t-},I_{t-},J_{t-}) - c(X_{t-},J_{t-})\,\one_{X_{t-} \in \partial E}.
 \end{equation*}
 We now  consider the inductive step: $1\leqslant k \in \N$,  and assume that $\P^{x,a_0,a_\Gamma}$-a.s.,
 \begin{align}
 	Y_{t}^{T,n,k} &= v^{T,n,k}(X_{t},I_{t}, J_{t})\label{Ynk}\\
 	Z_{t}^{T,n,k}(y,b,c)-\hat Z_{t}^{T,n,k}\,\one_{X_{t-} \in \partial E} &=v^{T,n,k}(y,b,c)-v^{T,n,k}(X_{t-},I_{t-},J_{t-}) - c(X_{t-},J_{t-})\,\one_{X_{t-} \in \partial E}.\label{Znk}
 \end{align}
 Then, plugging \eqref{Ynk}-\eqref{Znk} in \eqref{kBSDE} and computing the conditional expectation as before, 
 by the Markov property of $(X,I)$ we achieve that, $\P^{x,a_0,a_\Gamma}$-a.s.,
 $
 Y_{t}^{T,n,k+1}= v^{T,n,k+1}(X_{t},I_{t},J_t).
 $
 Then, applying the It\^o formula to $|Y_{t}^{T,n,k}-Y_{t}^{T,n}|^2$ and taking the supremum of $t$ between $0$ and $S$, one can show that
 $ \sperxa{\sup_{0 \leqslant t \leqslant S}\Big|Y_{t}^{T,n,k}-Y_{t}^{T,n}\Big|^2} \rightarrow 0$ as $k$ goes to infinity.
 Therefore, $v^{T,n,k}(x,a_0,a_\Gamma) \rightarrow v^{T,n}(x,a_0,a_\Gamma)$ as $k$ goes to infinity, for all $(x,a_0,a_\Gamma) \in E \times A_0 \times A_\Gamma$, from which it follows that, $\P^{x,a_0,a_\Gamma}$-a.s., 
 \begin{equation*}
 Y_{t}^{T,n,x,a_0,a_\Gamma}= v^{T,n}(X_{t},I_{t},J_t).
 \end{equation*}
This is the required identification property for $Y^{T,n,x,a_0,a_\Gamma}$. Letting $T\rightarrow\infty$, and then $n\rightarrow\infty$, we deduce property \eqref{Sec:PDP_ident_vdelta1} for $Y^{x,a_0,a_\Gamma}$.

\vspace{2mm}

\noindent\textbf{Step II.} \emph{The non-dependence of the function $v$ on its last arguments.} Notice that,
 by \eqref{Sec:PDP_Vstar_Y0} and \eqref{Sec:PDP_def_v}, $v$ coincides with the  value function $V^{\ast}$ of the randomized control problem introduced in \eqref{Sec:PDP_dual_value_function}.
 Therefore, to prove \eqref{Sec:PDP_vdelta_not_dep_a} we have  to show that $V^{\ast}(x,a_0,a_\Gamma)$ does not depend on $(a_0,a_\Gamma)$. In other words, given $(a_0,a'_0) \in A_0, (a_\Gamma,a'_\Gamma) \in A_\Gamma$, we have to prove that
\begin{equation}\label{V*a_a'}
V^*(x,a_0,a_\Gamma) \ = \ V^*(x,a_0',a_\Gamma').
\end{equation}
Notice that \eqref{V*a_a'} follows if we prove the following property of the cost functional: for every $\bm\nu=(\nu^0,\nu^\Gamma)\in\in \mathcal V$, there exist $(\nu^{0,\varepsilon}, \nu^{\Gamma,\varepsilon})_{\varepsilon}\in\mathcal V$ such that
\begin{equation}\label{Sec:PDP_claim_PDMP}
\lim_{\varepsilon \rightarrow 0^+} J(x,a'_0,a'_\Gamma,\nu^{0,\varepsilon},\nu^{\Gamma,\varepsilon}) = J(x,a_0,a_\Gamma,\nu^0,\nu^\Gamma).
\end{equation}
As a matter of fact, suppose that property \eqref{Sec:PDP_claim_PDMP} holds. Then, we deduce that
$V^{\ast}(x,a'_0,a'_\Gamma) \leq J(x,a_0,a_\Gamma,\nu^0,\nu^\Gamma)$,
 and by the arbitrariness of $(\nu^0,\nu^\Gamma)$, we conclude that
$V^{\ast}(x,a'_0,a'_\Gamma) \leq V^{\ast}(x,a_0,a_\Gamma)$,
from which we get \eqref{V*a_a'}. 

It remains to prove \eqref{Sec:PDP_claim_PDMP}. This can be done proceeding  as in the proof of Proposition 5.6 in \cite{BandiniPDMPsNoBordo}, that is as in the context of PDMPs with no jumps from the boundary, since the presence of predictable jumps does not induce here any additional  technical difficulty.

From now on, we suppose that the function $v$ is defined on $E$, as in \eqref{Sec:PDP_redefine_vdelta}. So, in particular, identity \eqref{Sec:PDP_ident_vdelta1} gives the representation formula \eqref{Sec:PDP_ident_vdelta}.

\vspace{2mm}

\noindent\textbf{Step III.} \emph{The function $v$ is bounded and continuous.} By   \eqref{Sec:PDP_Vstar_Y0}, \eqref{Sec:PDP_ident_vdelta}   and recalling the definition of $V^{\ast}$ in \eqref{Sec:PDP_dual_value_function},  
we have 
 $$
 v(x)=  V^\ast(x,a_0,a_\Gamma)
 = \inf_{\bm \nu \in \mathcal V} \spernuxa{\int_{(0,\,\infty)}e^{-\delta s} f(X_s,I_s)\, ds + \int_{(0,\,\infty)}e^{-\delta s} c(X_{s-},J_{s-})\, d p^\ast_s}.
 $$
  The boundedness of $v$ then  directly follows from the boundedness of $f$ and $c$. 
  In particular, $|v(x)|\leqslant \frac{M_f}{\delta} + C^\ast\,M_c$, for all $x\in E$.

 Let us now prove the continuity property of $v$. 
We proceed as in \cite{Da-Fa}, Section 5. Let $B(E)$ be the set of all bounded functions on $E$. Fix $(a_0, a_\Gamma) \in A_0 \times  A_\Gamma$, and  define the deterministic  operator $G: B(E) \rightarrow B(E)$ as $G \psi(x) :=  \inf_{\bm \nu \in \mathcal V} G_{\bm \nu} \psi(x)$, where 
$$
G_{\bm \nu} \psi(x) :=  \spernuxa{\int_{(0,\,T_1]}e^{-\delta s} f(X_s,I_s)\, ds + \int_{(0,\,T_1]}e^{-\delta s} c(X_{s-},J_{s-})\, d p^\ast_s + e^{-\delta T_1} \,\psi(X_{T_1})},
$$
with $T_1$   the first jump time of the PDMP $(X,I,J)$ under $\P^{x,a_0,a_\Gamma}_{\bm \nu}$.
 Set $t^\ast_{\bm \nu}(x) :=\inf \{t\geq 0 :\,X_t \in \partial E, (X_0,I_0,J_0)=(x,a_0,a_\Gamma),\,\,\P^{x,a_0,a_\Gamma}_{\bm \nu}{\textup{-a.s.}}\}$, and consider the sequence of Borel-measurable functions $(v_n)_{n \geqslant0}$ defined by 
\begin{align*}
&v_{n + 1}(x) = G\,v_n(x) := \inf_{\bm \nu\in \mathcal V} \left\{\int_0^{t^{\ast}_{\bm \nu}(x)}\chi^{\bm \nu}(s) f^{v_n}_0(X_s, I_s)\, ds + \chi^{\bm \nu}(t^{\ast}_{\bm \nu}(x)) F^{v_n}(X_{t^{\ast}_{\bm \nu}}, J_{t^{\ast}_{\bm \nu}})\right\},
\end{align*}
where $\chi^{\bm \nu}(s):= e^{-\delta s}\,e^{-\int_0^s\tilde \lambda^{\bm \nu}(t, X_t, I_t)\,dt}$ and, 
 for any $\psi \in B(E)$, 
\begin{align*}
f^{\psi}_0(X_s, I_s) &= f(X_s,I_s) + \int_E \psi(y)\,\lambda(X_s,I_s)\,Q(X_s,I_s,dy) \\
F^{\psi}(X_{s-}, J_{s-})&= c(X_{s-}, J_{s-}) + \int_E \psi(y)\, R(X_{s-}, J_{s-}, dy). 
\end{align*}
If we  prove  that $G$ is a two-stage contraction mapping, then by the strong Markov property of the PDMP $(X, I,J)$ it would follow that $v$ is the unique fixed point of $G$, and therefore $v(x) = \underset{n \rightarrow \infty}{\lim} v_n(x)$, see Corollary 5.6 in \cite{Da-Fa}. Then, the continuity property of $v$ in $E$ would follow from the existence of two monotone sequences of continuous functions converging to $v$, one from above and one from below, see Lemmas 5.9 and 5.10 in \cite{Da-Fa}.


It remains to prove that $G^2$ is a contraction in $E$. To this end,  it is enough to show that, for any $\psi_1, \psi_2 \in B(E)$,  $|G^2_{\bm \nu}\psi_1- G^2_{\bm \nu}\psi_2| \leq \rho ||\psi_1- \psi_2||$ for some constant $\rho < 1$, independent on $\bm \nu$, where  $||\psi|| = \max_{x \in E} \psi(x)$, $\psi \in B(E)$. 
Denoting by $T_2$ the second jump time of $(X,I,J)$, we have 
\begin{align*}
G^2_{\bm \nu} \psi(x) &:=  \spernuxa{\int_{(0,\,T_2]}e^{-\delta s} f(X_s,I_s)\, ds + \int_{(0,\,T_2]}e^{-\delta s} c(X_{s-},J_{s-})\, d p^\ast_s + e^{-\delta T_2} \,\psi(X_{T_2})},
\end{align*}
 so that $|G^2_{\bm \nu} \psi_1 - G^2_{\bm \nu} \psi_2| \leq \spernuxa{e^{-\delta T_2}}||\psi_1 - \psi_2||$. The fact that $\spernuxa{e^{-\delta T_2}} \leq \rho < 1$  is a consequence of  assumption \textbf{(H$0$)}, 
  see the proof of Proposition 46.17 in \cite{Da} for more details.
\endproof

We can finally state our main result. 

\begin{theorem}\label{Sec:PDP_THm_Feynman_Kac_HJB_2}
Let assumptions \textup{\textbf{(H$\textup{h$\lambda$QR}$)}},  {\bf (H$0$)}, {\bf (H$0'$)}, \textup{\textbf{(H$\lambda_0\lambda_\Gamma$)}} and \textup{\textbf{(H$\textup{fc}$)}} hold. Then, the function $v$ in \eqref{Sec:PDP_def_v} is a  viscosity  solution to \eqref{HJB_E}-\eqref{HJB_partialE}. Therefore, if assumptions \textup{\textbf{(HBB)}}, \textup{\textbf{(HBB')}} hold and $A_0$, $A_\Gamma$ are compact, by Theorem \ref{Thm_HJB_unique_viscosity_sol} we conclude that $v\equiv V$ and $V$ admits the Feynman-Kac representation formula
\begin{equation*} \label{Sec:PDP_Feynman-Kac}
V(x)=Y^{x,a_0,a_\Gamma}_0,\quad (x,a_0,a_\Gamma)\in E \times A_0 \times A_\Gamma.
\end{equation*}
\end{theorem}
Before proving Theorem \ref{Sec:PDP_THm_Feynman_Kac_HJB_2}, we recall the following useful technical result.
\begin{lemma}\label{L:equiv_def}
A function $u \in \mathbb C_b(\bar E)$ (resp. $w \in \mathbb C_b(\bar E)$)  is a sub- (resp. super-) solution to \eqref{HJB_E}-\eqref{HJB_partialE} if and only if, for any $\phi \in \mathbb C_b^1(\bar E)$, for any $x_0$ global maximum (resp. global minimum) point of $u- \phi$ (resp. $w- \phi$),
\begin{eqnarray*}
&&H^\phi(x_0,\phi(x_0),\nabla\phi(x_0)) \leqslant 0 \,\,\, \qquad \qquad \qquad \qquad \qquad \text{if} \,\, x_0 \in E,\\
&&\min\{H^\phi(x_0,\phi(x_0),\nabla\phi(x_0)),\phi(x_0) - F^\phi(x_0)\}  \leqslant 0 \quad \text{if} \,\, x_0 \in \partial E.
\end{eqnarray*}
\begin{eqnarray*}
\Big(\textup{resp.}
&&H^\phi(x_0,\phi(x_0),\nabla\phi(x_0)) \geqslant 0 \,\,\, \qquad \qquad \qquad \qquad \qquad \text{if} \,\, x_0 \in E,\\
&&\max\{H^\phi(x_0,\phi(x_0),\nabla\phi(x_0)),\phi(x_0) - F^\phi(x_0)\}  \geqslant 0 \quad \text{if} \,\, x_0 \in \partial E.\Big)
\end{eqnarray*}
\end{lemma}
\proof
See the proof of Proposition II.1 in \cite{Sa}.
\endproof

\noindent\emph{Proof (of Theorem \ref{Sec:PDP_THm_Feynman_Kac_HJB_2}).}\; Notice that, by Theorem \ref{Sec:PDP_THm_Feynman_Kac_HJB}, it is enough to check the viscosity sub- and super-solution properties for $v$ in the sense of Lemma \ref{L:equiv_def}. We split the proof into two steps.

\vspace{2mm}

\noindent\textbf{Viscosity subsolution property.} Let $\bar{x} \in \bar E$,  and let $\varphi \in C^1(\bar E)$ be a test function such that
 \begin{equation}\label{Sec:PDP_max_property}
 0= (v-\varphi)(\bar{x})=\max_{y \in \bar E} (v-\varphi)(y).
 \end{equation}
 \emph{Case 1: $\bar x \in E$}.
 Fix $(a_0,a_\Gamma) \in A_0 \times A_\Gamma$, set  $\eta=\frac{1}{2}\,d(\bar x, \partial E)$
 , and $ \tau :=\inf \{t \geqslant 0: |\phi(t,\bar x, a_0)-\bar x|\geqslant \eta\}$. 
Let $h>0$. 
Let $Y^{\bar x,a_0,a_\Gamma}$ be the unique maximal solution to \eqref{Sec:PDP_BSDE}-\eqref{Sec:PDP_BSDE_constraint1}-\eqref{Sec:PDP_BSDE_constraint2}  under $\P^{\bar x,a_0,a_\Gamma}$.
 We apply the It\^o formula to $e^{-\delta t} \,Y_t^{\bar x,a_0,a_\Gamma}$ between $0$ and $\theta := \tau \wedge h \wedge T_1$, where $T_1$ denotes the first jump time of $(X,I,J)$.
 From the constraints  \eqref{Sec:PDP_BSDE_constraint1}-\eqref{Sec:PDP_BSDE_constraint2} and the fact that $K$ is a nondecreasing process,   
  it follows that $\P^{\bar x,a_0,a_\Gamma}$-a.s.,
 \begin{align*}
 	Y^{\bar x,a_0,a_\Gamma}_0 &\leqslant e^{-\delta \theta_m}\,Y^{\bar x,a_0,a_\Gamma}_{\theta}
 	+ \int_{(0,\theta]} e^{-\delta r}\,f(X_{r}, I_r)\,dr
 	+ \int_{(0,\theta]} e^{-\delta r}\,c(X_{r-}, J_{r-})\,d p^\ast_r \\
 	&-\int_{(0,\theta]}e^{-\delta r}\,\int_{E \times A_0 \times A_\Gamma} Z^{\bar x,a_0,a_\Gamma}_r\,q(dr\,dy\,db\,dc).
 \end{align*}
 Applying the expectation with respect to $\P^{\bar x,a_0,a_\Gamma}$, from the identification property \eqref{Sec:PDP_ident_vdelta}, together with \eqref{Sec:PDP_max_property},
  it follows that 
 \begin{eqnarray*}
 	\varphi(\bar x)  \leqslant \E^{\bar x,a_0,a_\Gamma}\left[e^{-\delta \theta}\,\varphi(X_{\theta}) + \int_{(0,\theta]} e^{-\delta r}\,f(X_{r}, I_r)\,dr + \int_{(0,\theta]} e^{-\delta r}\,c(X_{r-}, J_{r-})\,d p^\ast_r\right].
 \end{eqnarray*}
 At this point, applying It\^o's formula to $e^{-\delta r}\,\varphi(X_{r} 
 )$ between $0$ and $\theta$, we get
 \begin{align}\label{Sec:PDP_sub_sol_int_ineq1}
 &  \frac{1}{h}\E^{\bar x,a_0,a_\Gamma}\left[\int_{(0,\theta]}e^{-\delta r}\,[\delta\,\varphi(X_{r}) -\mathcal L^{I_r} \varphi(X_{r})-f(X_{r}, I_r)]\,dr\right]\notag\\
 &\leq  \frac{1}{h}\E^{\bar x,a_0,a_\Gamma}\left[ e^{-\delta \theta}\,[\mathcal R^{J_{\theta-}} \varphi(X_{\theta-})+c(X_{\theta-}, J_{\theta-})]\,\one_{X_{\theta-} \in \partial E}\right],
 \end{align}
 where 
\begin{align}
\mathcal L^{I_r} \varphi (X_{r}) &:= h(X_r, I_r) \cdot \nabla \varphi(X_r) + \int_{E}(\varphi (y)-\varphi (X_{r}))\,\lambda(X_{r},I_r)\,Q(X_{r},I_r,dy),\label{LI}\\
\mathcal R^{J_{r-}} \varphi (X_{r-}) &:= \int_{E}(\varphi (y)-\varphi (X_{r-}))\,R(X_{r-},J_{r-},dy)\label{RJ}.
\end{align}
 Now we notice that, for every $r \in [0,\,\theta]$, $(X_{r-}, I_{r-}, J_{r-})= (\phi(r,\bar x, a_0),a_0,a_\Gamma)$, $\P^{\bar x, a_0,a_\Gamma}$-a.s., with  $\phi(r,\bar x, a_0) \in E$.  
 In particular the right-hand side of \eqref{Sec:PDP_sub_sol_int_ineq1} is zero. Taking into account the continuity on $E$ of the map   $z\mapsto \delta\,\varphi(z)-\mathcal L^{a_0}\varphi(z)-f(z,a_0)$, 
 we see that for any $\varepsilon>0$,
     \begin{align}\label{Sec:PDP_sub_sol_int_ineq2BIS}
 \frac{-\varepsilon + \delta\,\varphi(\bar x) -\mathcal L^{a_0} \varphi(\bar x)-f(\bar x,a_0)}{h} \,\E^{\bar x,a_0,a_\Gamma}\left[\frac{1-e^{-\delta \theta}}{\delta}\right]\leqslant 0.
 \end{align}
 Set $g(\theta):= \frac{1-e^{-\delta \theta}}{\delta}$, $\theta \in R_+$.
Recalling that the 
 the distribution density of $T_1$ under $\P^{\bar x,a_0,a_\Gamma}$ is given by  
 ($\tilde \lambda$ is the function introduced in \eqref{Sec:PDP_lambda_XI})
 $$
f_{T_1}(s)=  
\tilde \lambda(\phi(s,\bar x,a_0),a_0))\, \exp\left(-\int_0^s 
\tilde \lambda(\phi(r,\bar x,a_0),a_0)\,dr	\,\right)\,\one_{\phi(s,\bar x,a_0) \in E},
 $$
 we have
 \begin{align}
 &\frac{\E^{\bar x,a_0,a_\Gamma}\left[g(\theta)\right]}{h}= \frac{1}{h}\int_0^{h}
 g(s)\,f_{T_1}(s)\,ds + \frac{g(h)}{h}\,\P^{\bar x,a_0,a_\Gamma}[T_1> h]\nonumber\\
 &= \int_0^{h}
 \frac{1-e^{-\delta s}}{\delta\,h}\,\tilde\lambda(\phi(s,\bar x, a_0),a_0)\,e^{-\int_0^s \tilde \lambda(\phi(r,\bar x, a_0),a_0)\,dr}
 \,ds 
 + \frac{1-e^{-\delta h}}{\delta\,h}\,e^{-\int_0^{h}\tilde\lambda(\phi(r,\bar x, a_0),a_0)\,dr}.
 \label{Sec:PDP_2_terms}
 \end{align} 
 By the boundedness of $\lambda$, $\lambda_0$ and $\lambda_\Gamma$, it is easy to see that the two terms in the right-hand side of \eqref{Sec:PDP_2_terms} converge respectively to zero and one when $h$ goes to zero. 
 Thus, passing into the limit in \eqref{Sec:PDP_sub_sol_int_ineq2BIS} as $h$ goes to zero
 we obtain
\begin{equation*}
 \delta\,\varphi(\bar{x}) -h(\bar x , a_0) \cdot \nabla \varphi(\bar x) - \int_{E} (\varphi(y)- \varphi(\bar x)) \,\lambda(\bar x, a_0)\,Q(\bar x, a_0,dy)
 -f(\bar{x}, a_0) \leqslant 0.
 \end{equation*}
From the arbitrariness of  $a_0 \in A_0$, 
 we  conclude that $H^\varphi(\bar x, \varphi(\bar x), \nabla \varphi(\bar x)) \leq 0$.

\noindent \emph{Case 2: $\bar x \in \partial E$}.
If $\varphi(\bar x) - F^\varphi(\bar x)  \leq 0$ we have finished. Otherwise, suppose that $\varphi(\bar x) - F^\varphi(\bar x)  > 0$. 
We argue similarly to the  \emph{Case 1}. 
 Let $(x_m)_m$ in $E$ such that
 $
 x_m \underset{m \rightarrow \infty}{\longrightarrow} \bar{x}$.
 Fix $(a_0,a_\Gamma) \in A_0 \times A_\Gamma$.
 Let $\eta_m := \frac{1}{2} d(x_m, \partial E)$,  and 
$
 \tau_m :=\inf \{t \geqslant 0: |\phi(t,x_m, a_0)-x_m|\geqslant \eta_m\} 
$. 
Let $Y^{x_m,a_0,a_\Gamma}$ be the unique maximal solution to \eqref{Sec:PDP_BSDE}-\eqref{Sec:PDP_BSDE_constraint1}-\eqref{Sec:PDP_BSDE_constraint2}  under $\P^{x_m,a_0,a_\Gamma}$.
 We apply the It\^o formula to $e^{-\delta t} \,Y_t^{x_m,a_0,a_\Gamma}$ between $0$ and $\theta_m := \tau_m 
 \wedge T_1$, where $T_1$ denotes the first jump time of $(X,I,J)$ under $\P^{x_m,a_0,a_\Gamma}$.
 Proceeding as in \emph{Case 1}, we get 
 \begin{align}\label{Sec:PDP_sub_sol_int_ineq}
  &  \frac{1}{\tau_m} \E^{x_m,a_0,a_\Gamma}\left[\int_{(0,\theta_m]}\,e^{-\delta r}\,[\delta\,\varphi(X_{r}) -\mathcal L^{I_r} \varphi(X_{r})-f(X_{r}, I_r)]\,dr\right]\notag\\
 &\leq \frac{1}{\tau_m} \E^{x_m,a_0,a_\Gamma}\left[ \,e^{-\delta \theta_m}\,[\mathcal R^{J_{\theta_m-}} \varphi(X_{\theta_m-})+c(X_{\theta_m-}, J_{\theta_m-})]\,\one_{X_{\theta_m-} \in \partial E}\right],
 \end{align}
 where $\mathcal L^{I_r}$ and $\mathcal R^{J_{r-}}$ are the operators defined respectively in \eqref{LI} and \eqref{RJ}.
 Now we notice that, for every $r \in [0,\,\theta_m]$, $(X_{r-}, I_{r-}, J_{r-})= (\phi(r,x_m, a_0),a_0,a_\Gamma)$, $\P^{x_m, a_0,a_\Gamma}$-a.s., with  $\phi(r,x_m, a_0) \in E$. In particular the right-hand side of \eqref{Sec:PDP_sub_sol_int_ineq} is zero.
By the continuity of the map
  $\Gamma(z):=  \delta\,\varphi(z) -\mathcal L^{a_0} \varphi(z)-f(z, a_0)$,
 for any $\varepsilon>0$, there exists $l = l(\varepsilon)>0$ such that 
 $|\Gamma(y) - \Gamma(\bar x)| \leq \varepsilon$ if $|y-\bar x|\leq l(\varepsilon)$.
 Thus, for $\varepsilon$ fixed, let 
 $m = m(\varepsilon) \in \N$ such that, for any $m \geq m(\varepsilon)$, $\eta_m \leq \frac{1}{2} l(\varepsilon)$ and $|x_m - \bar x| \leq  \frac{1}{2} l(\varepsilon)$. By the triangle inequality, $|\phi(r,x_m, a_0)-\bar x| \leq l(\varepsilon)$. Therefore, for $m \geq m(\varepsilon)$, 
 \begin{align}\label{Sec:PDP_sub_sol_int_ineq2}
[-\varepsilon + \delta\,\varphi(\bar x) -\mathcal L^{a_0} \varphi(\bar x)-f(\bar x,a_0)] \,\E^{x_m,a_0,a_\Gamma}\left[\frac{1-e^{-\delta \theta_m}}{\delta \tau_m}\right]\leqslant 0.
 \end{align}
  Then, proceeding as in \emph{Case 1}, we get 
   \begin{align*}
 \E^{x_m,a_0,a_\Gamma}\left[\frac{1-e^{-\delta \theta_m}}{\delta \tau_m}\right] &= \int_0^{\tau_m}
 \frac{1-e^{-\delta s}}{\delta\,\tau_m}\,\tilde\lambda(\phi(s,x_m, a_0),a_0)\,e^{-\int_0^s \tilde \lambda(\phi(r,x_m, a_0),a_0)\,dr}
 \,ds \notag\\
 &+ \frac{1-e^{-\delta \tau_m}}{\delta\,\tau_m}\,e^{-\int_0^{\tau_m}\tilde\lambda(\phi(r, x_m, a_0),a_0)\,dr},
 \end{align*}
that goes to one as $m$ goes to infinity.
Thus, passing into the limit in \eqref{Sec:PDP_sub_sol_int_ineq2} as $m$ goes to infinity,
we  conclude also in this case  that $H^\varphi(\bar x, \varphi(\bar x), \nabla \varphi(\bar x)) \leq 0$ from the arbitrariness of  $a_0 \in A_0$.

\vspace{2mm}

\noindent\textbf{Viscosity supersolution property.} Let $\bar{x} \in \bar E$,  and let $\varphi \in C^1(\bar E)$ be a test function such that
 \begin{equation}\label{Sec:PDP_min_property}
 0= (v-\varphi)(\bar{x})=\min_{x \in \bar E} (v-\varphi)(x).
 \end{equation}
 \emph{Case 1: $\bar x \in E$}.
Notice that we can assume w.l.o.g. that $\bar{x}$ is a strict minimum of $v-\varphi$. 
As a matter of fact, one can subtract to $\varphi$ a positive cut-off function which behaves as $|x-\bar{x}|^2$ when $|x-\bar{x}|^2$ is  small, and that regularly converges to $1$ as  $|x-\bar{x}|^2$ increases to $1$.
  Then, for every $\eta >0$, 
 we can define
 \begin{equation}\label{Sec:PDP_min_beta}
 0 < \beta(\eta) := \inf_{x \in B^c(\bar{x},\eta) \cap \bar E} (v-\varphi)(x),
 \end{equation}
 where $B(\bar{x}, \eta):=\{y \in E: |\bar{x}-y|<\eta\}$.
 
 We will show the result by contradiction. Assume thus that
 $
 H^{\varphi}(\bar{x}, \varphi(\bar x), \nabla \varphi(\bar x)) < 0$.
 Then by the continuity of $H$, there exists $\eta >0$, $\beta(\eta) >0$ and $\varepsilon \in (0,\,\beta(\eta) \delta]$ such that
 $$
 H^{\varphi}(y, \varphi(y), \nabla \varphi(y)) \leqslant -\varepsilon,\quad \textup{for all}\,\, y \in B(\bar{x},\eta).
 $$
 Let us  fix $T >0$ and define $\theta := \tau \wedge T$, where $\tau =\inf \{t \geqslant 0: X_{t} \notin B(\bar{x}, \eta)\}$.
Moreover, let us fix  $(a_0,a_\Gamma) \in A_0 \times A_\Gamma$, and consider   the  solution  $Y^{n,\bar x,a_0,a_\Gamma}$ to the penalized \eqref{Sec:PDP_BSDE_penalized}, under the probability $\P^{\bar x,a_0,a_\Gamma}$.
 Notice that
 $\P^{\bar x,a_0,a_\Gamma}\{\tau=0\}= \P^{\bar x,a_0,a_\Gamma}\{X_0\notin B(\bar x, \eta)\}= 0$.
 
 We apply It\^o's formula to $e^{-\delta t} \,Y_t^{n,\bar x,a_0,a_\Gamma}$ between $0$ and 
 $\theta$.
 Then, proceeding as   in the proof of the representation formula \eqref{Sec:PDP_rep_Y_n}, we get the following inequality:
 \begin{equation}\label{Sec:PDP_ineq_Yn_nu}
 Y_{0}^{n,\bar x,a_0,a_\Gamma} \geqslant \inf_{\bm\nu \in \mathcal{V}^n} \E^{\bar x,a_0,a_\Gamma}_{\bm\nu}\left[e^{-\delta \theta 
 }\,Y_{\theta 
}^{n,\bar x,a_0,a_\Gamma} + \int_{(0,\theta]} 
 e^{-\delta r}\,f(X_{r}, I_r)\,dr +\int_{(0,\theta] 
} e^{-\delta r}\,c(X_{r-}, J_{r-})\,d p^\ast_r\right].
\end{equation}
Recall that $Y^{n,\bar x,a_0,a_\Gamma}$ converges decreasingly to the maximal solution $Y^{x_m,a_0,a_\Gamma}$ to the constrained BSDE \eqref{Sec:PDP_BSDE}-\eqref{Sec:PDP_BSDE_constraint1}-\eqref{Sec:PDP_BSDE_constraint2}. By the identification property \eqref{Sec:PDP_ident_vdelta}, together with  \eqref{Sec:PDP_min_property}
and \eqref{Sec:PDP_min_beta}, from  inequality \eqref{Sec:PDP_ineq_Yn_nu}  
we get that there exists a strictly  positive, predictable and bounded function $\bm\nu \in \mathcal{ V}$ such that
\begin{align*}
\varphi(\bar x) &\geqslant  \E^{\bar x,a_0,a_\Gamma}_{\bm\nu}\left[e^{-\delta \theta 
	}\,\varphi(X_{\theta 
}) + \beta\,e^{-\delta \theta 
}\,
\one_{\{\tau 
	\leqslant T\}}\right]\\
	&+\E^{\bar x,a_0,a_\Gamma}_{\bm\nu}\left[\int_{(0,\theta]} 
 e^{-\delta r}\,f(X_{r}, I_r)\,dr+\int_{(0,\theta]} 
 e^{-\delta r}\,c(X_{r-}, J_{r-})\,d p^\ast_r\right] - \frac{\varepsilon}{2\,\delta}.
\end{align*}
At this point, applying It\^o's formula to $e^{-\delta r}\,\varphi(X_{r})$ between $0$ and $\theta 
$, we get
\begin{align}
 & \E^{\bar x,a_0,a_\Gamma}_{\bm\nu}\left[\int_{(0,\theta]} 
 e^{-\delta r}\,[\delta\,\varphi(X_{r}) -\mathcal L^{I_r} \varphi(X_{r})-f(X_{r}, I_r)]\,dr\right]-\beta\,\E^{\bar x,a_0,a_\Gamma}_{\bm\nu}\left[e^{-\delta \theta 
}\,\one_{\{\tau 
	\leqslant T\}}\right]+ \frac{\varepsilon}{2}\notag\\
&-  \E^{\bar x,a_0,a_\Gamma}_{\bm\nu}\left[\int_{(0,\theta]} 
 e^{-\delta r}\,[\mathcal R^{J_{r-}} \varphi(X_{r-})+c(X_{r-}, J_{r-})]\one_{X_{r-} \in \partial E}\,dp^\ast_r\right]
 \geqslant 0,\label{Sec:PDP_super_sol_int_ineq}
\end{align}
where $\mathcal L^{I_r}$ and $\mathcal R^{J_{r-}}$ are defined respectively in \eqref{LI} and \eqref{RJ}.
Notice that, for $r \in [0,\,\theta] 
$, $X_{r-} \in B(\bar x, \eta) \subset E$. In particular, 
$[\mathcal R^{J_{r-}} \varphi(X_{r-})+c(X_{r-}, J_{r-})]\one_{X_{r-} \in \partial E} =0$. Moreover,
\begin{align*}
	\delta\,\varphi(X_{r}) -\mathcal L^{I_r} \varphi(X_{r})-f(X_{r}, I_r) &\leqslant \delta\,\varphi(X_{r}) -\inf_{b \in A_0}\{\mathcal L^{b} \varphi(X_{r})+f(X_{r},b)\}\\
	&= H^{\varphi}(X_{r}, \varphi(X_{r}), \nabla \varphi(X_{r}))\leqslant -\varepsilon,
\end{align*}
and therefore, from \eqref{Sec:PDP_super_sol_int_ineq} we obtain
\begin{align*}
	0 &\leqslant 
	 \frac{\varepsilon}{2\,\delta} +  \E^{\bar x,a_0,a_\Gamma}_{\bm\nu_m}\left[-\varepsilon \int_{(0,\theta]} 
	 e^{-\delta r}\,dr -\beta\,e^{-\delta \theta 
}\,
\one_{\{\tau 
	\leqslant T\}}\right] \\
&= - \frac{\varepsilon}{2\,\delta} +  \E^{\bar x,a_0,a_\Gamma}_{\bm\nu}\left[\left(\frac{\varepsilon}{\delta}- \beta\right)e^{-\delta \theta 
} \one_{\{\tau 
\leqslant T\}} + \frac{\varepsilon}{\delta} e^{-\delta \theta 
}\,\one_{\{\tau 
> T\}}\right]\\
&\leqslant 
- \frac{\varepsilon}{2\,\delta} +  \frac{\varepsilon}{\delta}\,\E^{\bar x,a_0,a_\Gamma}_{\bm\nu}\left[e^{-\delta T}\,\one_{\{\tau 
	> T\}}\right]
\leqslant   - \frac{\varepsilon}{2\,\delta} +  \frac{\varepsilon}{\delta}\,e^{-\delta T}.
\end{align*}
Letting $T$ go to infinity
we achieve the contradiction:
$0 \leqslant - \frac{\varepsilon}{2\,\delta}$.

\noindent \emph{Case 2: $\bar x \in \partial E$}.
As in the previous case, we can assume w.l.o.g. that $\bar{x}$ is a strict minimum of $v-\varphi$. 
  Then, for every $\eta >0$, 
 we can define
 \begin{equation*}
 0 < \beta(\eta) := \inf_{x \in \bar B^c(\bar{x},\eta) \cap \bar E} (v-\varphi)(x),
 \end{equation*}
 where $\bar B(\bar{x}, \eta):=\{y \in \bar E: |\bar{x}-y|<\eta\}$.

If $\varphi(\bar x) - F^\varphi(\bar x)\geq 0$ we have finished. Otherwise, assume that $\varphi(\bar x) - F^\varphi(\bar x)<0$.
We will show the result by contradiction. Assume thus that
 $
H^{\varphi}(\bar{x}, \varphi(\bar x), \nabla \varphi(\bar x)) < 0$. 
 Then by the continuity of $H$ and $F$, there exists $\eta >0$, $\beta(\eta) >0$ and $\varepsilon \in (0,\,\beta(\eta) \delta]$ such that
 $$
 H^{\varphi}(y, \varphi(y), \nabla \varphi(y)) \leqslant -\varepsilon,
 \quad \varphi(y) - F^\varphi(y) \leq -\varepsilon, \quad  \textup{for all}\,\, y \in \bar B(\bar{x},\eta).
 $$
 Let  us fix $T >0$ and define $\theta := \tau  \wedge T$, where $\tau =\inf \{t \geqslant 0: X_{t} \notin \bar B(\bar{x}, \eta)\}$.
 Arguing as in \emph{Case 1}, we get 
 \begin{align}
 & \E^{\bar x,a_0,a_\Gamma}_{\bm\nu}\left[\int_{(0,\theta]} 
 e^{-\delta r}\,[\delta\,\varphi(X_{r}) -\mathcal L^{I_r} \varphi(X_{r})-f(X_{r}, I_r)]\,dr\right]+ \frac{\varepsilon}{2}\label{Sec:PDP_super_sol_int_ineq2}\\
&-  \E^{\bar x,a_0,a_\Gamma}_{\bm\nu}\left[\int_{(0,\theta]} 
 e^{-\delta r}\,[\mathcal R^{J_r} \varphi(X_{r-})+c(X_{r-}, J_{r-})]\one_{X_{r-} \in \partial E}\,dp^\ast_r\right]-\beta\,\E^{\bar x,a_0,a_\Gamma}_{\bm\nu}\left[e^{-\delta \theta 
}\,\one_{\{\tau 
	\leqslant T\}}\right]\geqslant 0\notag
\end{align} 
for some $\bm \nu \in \mathcal V$.
Noticing that, for $r \in [0,\,\theta] 
$, 
\begin{align*}
	\delta\,\varphi(X_{r}) -\mathcal L^{I_r} \varphi(X_{r})-f(X_{r}, I_r) &\leqslant \delta\,\varphi(X_{r}) -\inf_{b \in A_0}\{\mathcal L^{b} \varphi(X_{r})+f(X_{r},b)\}
	\leqslant -\varepsilon,
\\
-\left(\mathcal R^{J_r} \varphi(X_{r-})+c(X_{r-}, J_{r-}) \right)
&\leq -\min_{d \in A_\Gamma}\{\mathcal R^{d} \varphi(X_{r-})+c(X_{r-}, d)
\leq -\varepsilon,
\end{align*}
from \eqref{Sec:PDP_super_sol_int_ineq2}
we obtain 
\begin{align*}
	0 &\leqslant   \frac{\varepsilon}{2\,\delta} +  \E^{\bar x,a_0,a_\Gamma}_{\bm\nu}\left[-\varepsilon \int_{(0,\theta]} 
	 e^{-\delta r}\,dr-\varepsilon \int_{(0,\theta]} 
	 e^{-\delta r}\,dp^\ast_r -\beta\,e^{-\delta \theta 
}\,
\one_{\{\tau 
	\leqslant T\}}\right] \\
&= - \frac{\varepsilon}{2\,\delta} +  \E^{\bar x,a_0,a_\Gamma}_{\bm\nu}\left[\left(\frac{\varepsilon}{\delta}- \beta\right)e^{-\delta \theta 
} \one_{\{\tau 
\leqslant T\}} + \frac{\varepsilon}{\delta} e^{-\delta \theta 
}\,\one_{\{\tau 
> T\}}\right] -\varepsilon \,\E^{\bar x,a_0,a_\Gamma}_{\bm\nu}\left[ \int_{(0,\theta]} 
	e^{-\delta r}\,dp^\ast_r\right]\\
&\leqslant
 - \frac{\varepsilon}{2\,\delta} +  \frac{\varepsilon}{\delta}\,\E^{\bar x,a_0,a_\Gamma}_{\bm\nu}\left[e^{-\delta T}\,\one_{\{\tau 
	> T\}}\right]\leqslant  - \frac{\varepsilon}{2\,\delta} +  \frac{\varepsilon}{\delta}\,e^{-\delta T}.
\end{align*}
 Letting $T$ go to infinity we get the contradiction: $0 \leq -\frac{\varepsilon}{2\,\delta}$.
\qed

\appendix 
\renewcommand\thesection{Appendix} 
\section{} 
\renewcommand\thesection{\Alph{subsection}} 
\renewcommand\thesubsection{\Alph{subsection}}

\subsection{Proof of Proposition \ref{P_Prob_infty}}
\label{S:AppP}

By the Radon-Nikodym theorem, there exist three nonnegative functions $d_1$, $d_2$, $d_3$ defined on $\Omega \times [0,\,\infty) \times E \times A_0 \times A_\Gamma$, $\mathcal P \otimes \mathcal E \otimes \mathcal A_0 \otimes \mathcal A_\Gamma$, such that
\[
 d \tilde{p}^{\bm\nu} = (\nu^0 \, d_1 + \nu^\Gamma \, d_2  + d_3)\, d\tilde{p},
\]
with $d_1+d_2+d_3=1$, $\tilde p$-a.e., and
 \begin{align}
 &d_1(t,y,b,c)\, \tilde{p}(dt\,dy\,db\,dc) =  \lambda_0(db) \, \delta_{\{ X_{t-} \}}(dy)\, \delta_{\{ J_{t-} \}}(dc) \,	\one_{X_{t-} \in E}\, dt,\label{d1}
 	\\
 &d_2(t,y,b,c)\, \tilde{p}(dt\,dy\,db\,dc) =  	\lambda_\Gamma(dc) \, \delta_{\{ X_{t-} \}}(dy) \, \delta_{\{ I_{t-} \}}(db)\, \one_{X_{t-} \in E}\,dt,\label{d2}
 	\\
 &d_3(t,y,b,c)\, \tilde{p}(dt\,dy\,db\,dc) =
 	 \lambda(X_{t-},\,I_{t-})\,Q(X_{t-},\,I_{t-},\,dy) \, \delta_{\{ I_{t-} \}}(db) \, \delta_{\{ J_{t-} \}}(dc)\, \one_{X_{t-} \in E}\,dt + \notag\\
 	 &\qquad \qquad \qquad \qquad \qquad \quad\;\; + R(X_{t-},\,J_{t-},\,dy) \, \delta_{\{ I_{t-} \}}(db) \, \delta_{\{ J_{t-} \}}(dc)\, \one_{X_{t-} \in \partial E}\,d p^\ast_t.\notag
 \end{align}

 \begin{Remark}\label{R_d1d2d3}
Notice that, by construction, $d_1(t,y,b,c) \,\one_{X_{t-} \in \partial E} =d_2(t,y,b,c) \,\one_{X_{t-} \in \partial E} =0$, and $d_3(t,y,b,c) \,\one_{X_{t-} \in \partial E} =\one_{X_{s-} \in \partial E}$.
\end{Remark}

For any $\bm\nu\in \mathcal{V}$, define the Dol\'eans-Dade exponential local martingale
 $L^{\bm\nu}$ defined by 
 \begin{align*}
 L_t^{\bm\nu} &= 
  e^{\int_0^t\int_{A_0}(1 - \nu^0_r(b) )\lambda_0(db)\,dr} \,e^{\int_0^t\int_{A_\Gamma}(1 - \nu^\Gamma_r(c) )\lambda_\Gamma(dc)\,dr} \cdot\\
 &\cdot\prod_{n \geqslant 1: T_{n} \leqslant t}(\nu^0_{T_{n}}(A^0_n)\,d_1(T_{n},E_n,A^0_n, A^\Gamma_n)+\nu^\Gamma_{T_{n}}(A_n^\Gamma)\,d_2(T_{n},E_n,A^0_n, A^\Gamma_n) + d_3(T_{n},E_n,A^0_n, A^\Gamma_n)),\notag
 \end{align*}
 for all $t\geq 0$. Notice that, when $(L^{\bm\nu}_t)_{t \geq 0}$ is a true martingale, for every time $T>0$ we can define  a probability measure $\P^{x,a_0,a_\Gamma}_{\bm\nu, T}$ equivalent to $\P^{x,a_0,a_\Gamma}$ on $(\Omega,\,\mathcal{F}_T)$ by
 \begin{equation*}
 \P^{x,a_0,a_\Gamma}_{\bm\nu, T}(d\omega)=L_T^{\bm\nu}(\omega)
 \,\P^{x,a_0,a_\Gamma}(d\omega).
 \end{equation*}
\begin{lemma}\label{P:JacodCond}
 When $(L^{\bm\nu}_t)_{t \geq 0}$ is a true martingale, for every $T>0$, the restriction of the random measure $p$ to $(0,T]\times E\times A_0\times A_\Gamma$
 admits $\tilde{p}^{\bm\nu}=(\nu^0\,d_1+ \nu^\Gamma\,d_2 + d_3)\,\tilde{p}$
 as compensator   under $\P^{x,a}_{\bm\nu, T}$.
 \end{lemma}
 \proof
 We shall prove that 
\begin{equation}\label{cond_T4.5_J}
	 \hat {\bar \nu}_t=1  \,\,\, \textup{whenever}  \,\,\,\alpha_t=1,
\end{equation}
with $ \alpha_t := \tilde{p}(\{t\}\times E \times A_0 \times A_{\Gamma})$, and 
 \begin{align*}
 &\bar \nu_t(y,b,c) :=\nu^0_t(b)\,d_1(t,y,b,c)+ \nu^\Gamma_t(c)\,d_2(t,y,b,c) + d_3(t,y,b,c), \\
  &\hat{\bar \nu}_t:=\int_{E \times A_0 \times A_{\Gamma}}\bar \nu_t(y,b,c)\,\tilde{p}(\{t\}\times dy\,db\,dc).
 \end{align*}
Indeed, if  condition \eqref{cond_T4.5_J} holds, then the result would be a direct application of  Theorem 4.5  \cite{J}.
Let us thus show the validity of \eqref{cond_T4.5_J}. 
To this end, we start by noticing that, by Remark \ref{R_d1d2d3}, 
\begin{equation}\label{nuind}
\bar \nu_s(y,b,c)\,\one_{X_{s-} \in \partial E}= \one_{X_{s-} \in \partial E}.
\end{equation}
Moreover, \eqref{tildep_t} implies
\begin{align*}
\int_{E}\bar \nu(t,y,b,c)\,\tilde{p}(\{t\}\times dy\,db\,dc) = \int_{E}\bar \nu(t,y,I_{s-},J_{s-})\,R(X_{s-},J_{s-},dy)\,\one_{X_{s-} \in \partial E}=\one_{X_{s-} \in \partial E},
 \end{align*}
 where the latter equality follows from \eqref{nuind}.
 On the other hand, by \eqref{DeltaA} we have $\alpha_t=\one_{X_{s-} \in \partial E}$, and condition \eqref{cond_T4.5_J} follows.
 \endproof
We can now state the main result of this section.
\begin{proposition}\label{Sec:PDP_lemma_P_nu_martingale}
 	Let assumptions {\bf (Hh$\lambda$QR)}
 	and  {\bf (H$\lambda_0\lambda_\Gamma$)} hold. Then, for every  $(x,a_0,a_\Gamma)\in E\times A_0 \times A_\Gamma$  and
 	$\bm\nu \in \mathcal{V}$, under $\P^{x,a_0, a_\Gamma}$ the process $(L^{\bm\nu}_t)_{t \geq 0}$ is a   martingale. Moreover, for any $T>0$,  $L_T^{\bm\nu}$ is square integrable,
 	and, for every   $H \in \mathcal{G}_{\textbf{x},\textbf{a}_0, \textbf{a}_\Gamma}^{\textbf{2}}(\textup{q};\textup{\textbf{0,\,T}})$, 
 	the process 
 \begin{equation}\label{Mnu}
   M^{\bm\nu}_t:=\int_{(0,t]}\int_{E\times A_0 \times A_{\Gamma}} H_s(y,b,c)\,q^{\bm \nu}(ds\,dy\,db\,dc), \quad t\in[0,T],
\end{equation}
is a square integrable $\P^{x,a_0,a_\Gamma}_{\bm\nu, T}$-martingale on $[0,T]$. Finally, there exists a unique probability $\P^{x,a_0,a_\Gamma}_{\bm\nu}$ on $(\Omega, \mathcal F_{\infty})$, under which $\tilde{p}^{\bm\nu}$ in \eqref{Sec:PDP_dual_comp} is the compensator of $p$ in \eqref{Sec:PDP_p_dual} on $(0,\,\infty)\times E \times A_0 \times A_\Gamma$, and such that, for any $T >0$, the restriction of $\P^{x,a_0,a_\Gamma}_{\bm\nu}$ on $(\Omega, \mathcal F_T)$ is $\P^{x,a_0,a_\Gamma}_{\bm\nu,T}$.
\end{proposition}
\proof The first part of the result is a consequence of Theorem 5.2 in \cite{J}. 
 The square integrability property of $L^{\bm \nu}$ can be proved arguing as in  the proof of Lemma 3.2 in \cite{BandiniFuhrman}.
Moreover, Proposition 3.71-(a) in \cite{JB} implies that 
the stochastic integral $\int_{(t,T]} \int_{E\times A_0 \times A_{\Gamma}} H_s(y,b,c) \, q(ds\,dy\,db\,dc)$
   is well-defined, and, by  Proposition 3.66 in \cite{JB}, the process 
   $
   M_t:=\int_{(0,t]}\int_{E\times A_0 \times A_{\Gamma}} H_s(y,b,c)\,q(ds\,dy\,db\,dc)$, $t\in[0,T]$,
   is a square integrable  $\P^{x,a_0,a_\Gamma}_T$-martingale. Using the Burkholder-Davis-Gundy and
 Cauchy-Schwarz inequalities, together with the square integrability of  $L_T^\nu$, we see that \eqref{Mnu} is a square integrable  $\P^{x,a_0,a_\Gamma}_{T,\bm \nu}$-martingale. Finally, the probability measure $\P^{x,a_0,a_\Gamma}_{\bm\nu}$ on $(\Omega, \mathcal F_{\infty})$ can be constructed as usual by means of the Kolmogorov extension theorem (see e.g. Theorem 1.1.10  in \cite{Stroock_Varadhan}).
\endproof

\subsection{Some useful properties of the space $\mathcal{G}^{\textbf{2}}_{\textbf{x},\textbf{a}_0, \textbf{a}_\Gamma}(\textup{q};\textbf{0,\,T})$}
\label{S:AppG}

Following \cite{JB}, we define the random sets:
   \begin{align}
     D&:=\{(\omega, t) : \, p (\omega, \{t\}\times E \times A_0 \times A_\Gamma)>0 \},\label{D}\\
    J&:=\{(\omega, t) : \, \tilde p (\omega, \{t\}\times E \times A_0 \times A_\Gamma)>0 \},\notag\\
   K&:=\{(\omega, t) : \, \tilde p (\omega, \{t\}\times E \times A_0 \times A_\Gamma)=1 \}.\label{K}
    \end{align}
Notice that, by \eqref{DeltaA},   for any $t \geq 0$, 
\begin{equation}\label{J=K}
J = K =\{(\omega,t): \Delta A_t(\omega) = 1\}=\{(\omega,t): \,X_{t-}(\omega) \in \partial E\}.
\end{equation}

    	\begin{lemma}\label{L_simplifiedNorm} 
    	Let $H \in \mathcal{G}^{\textbf{2}}_{\textbf{x},\textbf{a}_0, \textbf{a}_\Gamma}(\textup{q};\textup{\textbf{0}},\,\textup{\textbf{T}})$. Then
    	\begin{itemize}
    	\item[(i)] $||H||^2_{\mathcal{G}_{\textbf{x},\textbf{a}_0, \textbf{a}_\Gamma}^{\textbf{2}}(\textup{q};\textup{\textbf{0,\,T}})} =$
\begin{align}\label{goodnorm}
 	& \sperxa{\int_{(0,\,T]}\int_{E\times A_0 \times A_\Gamma} \big|H_t(y,b,c) - \hat H_t\big|^2\,\tilde p(dt\,dy\,db\,dc) + \sum_{s \in (0,\,T]} |\hat H_s|^2 (1-\Delta A_t)	},
	\end{align}
	with $\hat H$ as in \eqref{Zhat}.
   \item[(ii)]$||H||^2_{\mathcal{G}_{\textbf{x},\textbf{a}_0, \textbf{a}_\Gamma}^{\textbf{2}}(\textup{q};\textup{\textbf{0,\,T}})}=\sperxa{\sum_{s \in (0,\, T]}
\bigg|\int_{E \times A_0 \times A_\Gamma}\,H_{s}(y,b,c)\,q(\{s\}\times  dy\,db\, dc)\bigg|^2}$.
		\end{itemize}
   \end{lemma}
   \proof Regarding item (i), the right-hand side of  \eqref{goodnorm} can be rewritten as  
    \begin{align*}
 	\sperxa{\int_{(0,\,T]}\int_{E \times A_0 \times A_\Gamma} \big|H_t(y,b,c) - \hat H_t\,\one_{J}(t)\big|^2\,\tilde p(dt\,dy\,db\,dc)
		 \ + \sum_{0 < t \leq T} \big|\hat H_t\big|^2\big(1 - \Delta A_t\big)\,\one_{J\setminus K}(t)},
		\end{align*}
	see e.g. Remark 2.6 in \cite{BandiniRusso2}, 
     and we  conclude by  \eqref{J=K}. Concerning item (ii), we have 
\begin{align*}
&
\bigg|\int_{E \times A_0 \times A_\Gamma}\,H_{s}(y,b,c)\,q(\{s\}\times  dy\,db\, dc)\bigg|^2 = \int_{E \times A_0 \times A_\Gamma}\,|H_{s}(y,b,c)- \hat H_s|^2\,p(\{s\}\times  dy\,db\, dc)\\
&- \int_{E \times A_0 \times A_\Gamma}\,|\hat H_s|^2\,p(\{s\}\times  dy\,db\, dc)+ |\hat H_s|^2 + 2 \int_{E \times A_0 \times A_\Gamma}\,\hat H_s\,H_s(y,b,c)\,p(\{s\}\times  dy\,db\, dc)\\
&- 2\,\hat H_s\int_{E \times A_0 \times A_\Gamma}\,H_s(y,b,c)\,p(\{s\}\times  dy\,db\, dc), 
\end{align*}
and the conclusion follows from item (i),  \eqref{DeltaA} and \eqref{Zhat}.
   \endproof	
For any stopping time $\tau$, we denote by  $[[\tau]]$ the random set $\{(\omega, \tau(\omega))\}\subset \Omega \times [0,\,\infty]$.
    \begin{lemma}\label{P_simplifiedNorm}
    Let $D$ be the random set in \eqref{D}.
Then  $D = K \cup (\cup_n[[T_n^i ]])$ up to an evanescent set, where $(T_n^i)_n$ are totally inaccessible times such that $[[ T_n^i ]] \cap [[ T_m^i ]] = \emptyset$, $n \neq m$.
 \end{lemma}
 \proof
Set $p^c := p \one_{J^c}$ and $\tilde p^c := \tilde p \one_{J^c}$. The measure $\tilde p^c$ is the compensator of $p^c$, see paragraph b) in \cite{JB}. We have
$$
\tilde p^c(ds \,dy\,db\,dc) = \tilde Q(X_{s−},I_{s-}, J_{s-},dy\,db\,dc)\,\one_{X_{s-} \in E}\, ds.
$$
We remark that
 $D \cap J^c = \{(\omega,t) : p^c(\omega,\{t\}\times \R) > 0\}$.
On the other hand, being $\tilde p^c$ absolutely continuous with respect to the Lebesgue measure, we have $\{(\omega,t) : p^c(\omega,\{t\}\times \R) > 0\}= \cup_n[[T_n^i]]$, for some $(T_n^i)_n$ totally inaccessible times, see, e.g., Assumption (A) in \cite{CFJ}. Therefore, since  $J = K$ by \eqref{J=K}, we have $D = K \cup (D \cap J^c)$, and the conclusion follows.
\endproof

\small

\paragraph{Acknowledgements.}
The  author  benefited
 from the support of the  Italian MIUR-PRIN 2015-16 "Deterministic and stochastic evolution equations".

\end{document}